\documentclass[preprint]{elsarticle}
\usepackage[colorlinks]{hyperref}
\usepackage{lineno,hyperref}
\usepackage{amsmath}
\usepackage{amsthm}
\usepackage{mathtools}
\usepackage{lipsum}
\usepackage{color}
\usepackage{ulem}
\usepackage{bbm}
\usepackage{graphicx}
\usepackage[utf8]{inputenc}
\usepackage{todonotes}
\newtheorem{theorem}{Theorem}
\newtheorem{definition}{Definition}
\newtheorem{observation}{Observation}
\newtheorem{assumption}{Assumption}
\modulolinenumbers[5]

\journal{}

\begin{document}

\begin{frontmatter}

\title{Discriminatory Price Mechanism for Smart Grid}

\author{Diptangshu Sen, Kushaagra Goyal, Varun Ramamohan, Arnob Ghosh\corref{mycorrespondingauthor}}
\address{Department of Mechanical Engineering, Indian Institute of Technology Delhi\\
Hauz Khas, New Delhi 110016, INDIA}
\cortext[mycorrespondingauthor]{Corresponding author. Email id: arnob.ghosh@imperial.ac.uk}




\begin{abstract}
We consider a scenario where a retailer can set different prices for different consumers in a smart grid. The retailer's objective is to maximize the revenue, minimize the operating cost, and maximize the consumer's welfare. The retailer wants to optimize a convex combination of the above objectives using price
signals specific to each consumer. However, variability in unit prices across consumers is bounded by a parameter $\eta$, hence limiting the discrimination. We formulate the pricing problem as a Stackelberg game where the retailer is the leader and consumers are followers. Since the retailer's optimization problem turns out to be non-convex, we convexify it via relaxations. We provide performance guarantees for the relaxations in the asymptotic sense (when number of consumers tends to $\infty$). Further, we show that despite the variability in pricing, the pricing scheme proposed by our model is \textit{fair} as higher prices are charged to consumers who have higher willingness for demand. We extend our analysis to the scenario where consumers can feed energy back to the grid via net-metering. We show that our pricing policy promotes \textit{fairness} even in this scenario as prosumers who contribute more to the grid, are given large cuts on buying rates. The policy is also found to incentivize more prosumers to invest in renewable energy, thus encouraging sustainability. 
\end{abstract}

\begin{keyword}
OR in energy, Pricing, Stackelberg games, Smart grid
\end{keyword}

\end{frontmatter}


\section{Introduction}
\subsection{Motivation} 
The advent of home automation has enabled consumers to control their consumption depending on the prices. Users' load can be categorized in two parts: elastic, and inelastic. Elastic load refers to the energy consumed by appliances whose usage can be deferred or optimized depending on the price, i.e, loads which are price-elastic. Charging of electric vehicles is a typical example because the user can optimally choose the charging schedule depending on the price points across the day. HVAC systems also constitute elastic loads because the user is directly influenced by the price signal while choosing the setpoint temperature. On the other hand, inelastic load refers to the load contributed by essential appliances like lights, induction-based cooking ovens etc. which often cannot be foregone and hence are price-inelastic. 

Using the advanced metering infrastructure installed at the users' premises, the retailer or utility company\footnote{We use the terms retailer and utility company interchangeably.} now can exercise better control on the total elastic load by varying the prices. For example, when the demand is high, the retailer can set high prices in order to dissuade the users from consuming high amounts of energy.  However, in the absence of a proper pricing mechanism, the full potential of the elastic load cannot be utilized. In recent times, demand response mechanisms have been proposed where users are compensated for reducing their consumption during periods of peak load \cite{liu2014pricing}. However, such price mechanisms have not reached their full potential because the compensation to the users is often not lucrative enough to induce changes in usage patterns. 

One of the primary reasons of failure of the demand-response price mechanism is that the price is same for all users. Different users have different valuations for consumption which means that different levels of incentivization are required to persuade them to participate in the demand-response programs. Further, researchers have argued that charging different prices to different users can in fact increase the efficiency\cite{sean}. It is also a practically implementable solution now due to the smart metering infrastructure as the retailer can directly communicate different prices to different users over the smart meter. Though such a pricing mechanism can be termed as `discriminatory', we show in this paper that it can actually promote `fairness' in terms of energy distribution in the community as the retailer is now charging higher prices to users who have higher willingness for demand and lower prices to those with low willingness for demand.\footnote{User $A$ is said to have a higher willingness for demand than another user $B$ when $A$ would consume more energy than $B$ for the same price per unit.} 

In the scenario where renewable energies are integrated with the power grid, discriminatory pricing is even more important. It is known that with higher penetration of renewable energy, there is increased variability in power supply and therefore, robust demand-response programs are required to match demand with supply at all times. Different users need to be incentivized at different rates so that they are keen to contribute to the main grid. Thus, discriminatory pricing can help to realize the enormous potential of renewable energy. 


\subsection{Our Approach} We seek to answer whether allowing prices to vary within a certain limit across different consumers can result in gains in social welfare and/or retailer's profit. We employ a game-theoretic framework to explore the same. We start with the scenario where users do not have access to renewable energy and are dependent entirely on the grid for meeting their demands. We formulate the problem as a Stackelberg game where the retailer sets the price and each consumer chooses how much to consume in response to the price. The retailer optimizes a weighted sum of her profit and consumers' welfare.  We show that the optimization problem of the retailer is non-convex and we address this by developing convex relaxations of the same. We then derive the conditions under which the relaxations give the same price as the original formulation. We also provide bounds for the optimality gap between the different relaxations and the original formulation in the asymptotic sense (when the number of users in the community tends to $\infty$).

In our formulation, we consider that the retailer is allowed to discriminate among consumers by charging them at different rates. However, we restrict the discrimination by an amount $\eta$. Even though we allow prices to be different across consumers, {\em we show that the prices achieved in both relaxations are fair}, i.e., the retailer charges higher prices only to those who have higher willingness for demand (Theorem \ref{thm:fair}). We also empirically evaluate how variation in $\eta$ can impact the retailer profits and consumer payoffs. 

Subsequently, we consider the scenario where the user generates renewable energy in-house and can feed some of it back to the grid. We investigate the net-metering price mechanism where the retail rate and the sell-back rate (rate at which users are compensated for selling back to the grid) are the same. Thus, if the retailer sets a higher price, 'prosumers' (consumers who can produce) can be incentivized to sell back more. High sell-back will reduce the operating costs for the retailer but their revenues will also decrease substantially. Hence it is not \textit{a priori} clear which price will maximize the retailer's objective. We formulate this problem of determining the optimal net-metering price for the retailer as an optimization problem and again develop convex relaxations for the same. In addition to being \textit{fair}, we also show that discriminatory pricing has the potential to incentivize prosumers to invest in renewable energy. 

\subsection{Literature Review} Demand response mechanisms using distributed optimization have been proposed in multiple settings for better demand-supply matching \cite{li,chang,chang2,gadh}. \cite{li} is of particular relevance because it proposes an optimal dynamic pricing scheme for demand response, similar to us. The prices are computed iteratively by the utility company by obtaining revised consumption estimates and load schedules from individual consumers in each step, until the price converges to the optimum. Since a lot of information needs to be exchanged for a price computation, it is vulnerable from a data security perspective. So need arises for a  price mechanism which has a low communication overhead. Alternatively, real time pricing mechanisms have been proposed in \cite{ma2014distributed,liu2014pricing,tsui,samadi2010optimal,wang} where the retailer sets the prices in a dynamic setting by estimating the total consumption. However, real-time pricing is anticipatory in nature for the consumers because they need to anticipate the prices while optimizing their consumption. Uncertainty in prices often leads to instabilities in the users' response. \cite{tsitsiklis2015pricing} is a unique work where the ancillary costs associated with demand variability are considered in the analysis. Additionally, they show that pricing at the marginal rate of generation does not maximize social welfare. However it is to be noted that discriminatory pricing or pricing for net-metering scenarios have not been considered in the scope of any of the above papers. 

Since all stakeholders are strategic and consumers are often not cooperating with each other, game theoretic models are best suited for modeling the interactions in retail electricity markets. There is usually a hierarchy in the interaction because the retailer always moves first in the game, so the Stackelberg game framework is commonly used in literature. Stackelberg games where the retailer and the consumers constitute the leader-follower pair,  have been considered in \cite{wei2014energy,chen2011innovative,saad2012game,maharjan,yu}. Other kinds of game theory models like cake-cutting games have also been explored elsewhere in the pricing setting (\cite{tushar2016price}). 

In recent times, due to the rapid integration of renewable energy sources into the demand side of the power grid, there is a need to redesign pricing schemes. \cite{zafar2018prosumer,liuprosumer,kanchev2011energy} discuss energy management systems for microgrids involving users who can contribute energy to the grid (also known as 'prosumers'). However, these papers consider different buying and selling prices which significantly increases the communication overhead. In comparison, we consider a pricing mechanism where the retail and sell-back rates are the same. This mechanism, popularly known in literature as \textit{net-metering} is easier to implement because it is one-shot and repeated exchange of information is not required. Net-metering has been studied extensively in literature and their benefits have been demonstrated in \cite{sajjad2015net, poullikkas2013comparative}. However, the above papers are aimed primarily at PV-residential systems while our setting is much more general and focuses on discriminatory pricing within the purview of net-metering.

The advent of smart metering has enabled the retailer to access the real-time consumption of each user and also set different prices for different users for better efficiency. There are several papers in literature which explore the concept of discriminatory pricing in the power grid. It has been shown by \cite{sean} that discriminatory pricing has the potential to increase efficiency. \cite{tushar2016price} also investigates discrimination in prices, but in a completely different setting from ours, where different prosumers in a smart community choose different prices at which they wish to sell energy to the SFC (shared facility controller). The problem is modelled using a cake-cutting game and the pricing scheme is shown to result in an envy-free marketplace. \cite{tushar} also considers a discriminatory pricing regime but the optimal prices are solved in a distributed manner. Our work differs from it because our pricing mechanism is one-shot. Further, unlike any of the above papers, we have shown that our mechanism is \textit{fair} (because it charges higher prices to users who have higher willingness for demand) and the fairness property is preserved even in the net-metering scenario. We also empirically study the impact of the discrimination limiting parameter $\eta$ on the retailer's revenues and the user's welfare. 

\subsection{Original Contributions} The main contributions of our work are the following:
\begin{itemize}
    \item We consider a piecewise-linear price response function for the consumer to capture the notion of inelastic and elastic demand. This makes our pricing model more realistic compared to the demand response model which studies optimal compensation for reducing the elastic load only. As a result, our formulation becomes non-convex, on the other hand all the state-of-the art papers consider a convex optimization problem formulation. We relax the problem to make it convex and provide performance guarantees for relaxed problems in the asymptotic sense.
    \item We develop a stylized model where the retailer's objective is to maximize the consumer's welfare in addition to maximizing her revenue and minimizing operation costs. We propose a discriminatory pricing scheme where the retailer can charge different users at different rates. We show that, in spite of being discriminatory, our pricing scheme is `fair' as consumers who have higher willingness for demand, are priced higher. It also leads to a `fairer' distribution of energy in the community in the sense that the total elastic load on the grid remains practically the same, but the variance in consumption across users is reduced. We empirically evaluate the impact of the level of discrimination on the retailer's objectives.
    \item We extend our analysis to the scenario where consumers may have access to renewable resources and can feed back energy to the grid. We formulate prices for a net metering scenario where the retail rate and sell-back rate are identical. Discriminatory pricing is found to be beneficial in this setting as well. The retailer is forced to offer higher rates to reluctant prosumers to induce sell-back, while prosumers who contribute high amounts to the grid, are given higher discounts on retail rates.
\end{itemize} 

\section{System Model}
\subsection{Entities \& Interactions}
There is a retailer or utility company $R$ that purchases electricity from the wholesale  market and supplies to a community of $N$ consumers. Time is discretized into $K$ periods. In each of the $k^{th}$ periods ($k \in \{1,2,...,K$), the retailer communicates their prices to the consumers based on which they choose their demands for the said period. Note that the duration of a period can be of any magnitude; however, the computational and communication overheads for the retailer increase as the duration of each time period decreases. 

Total demand for a household in a period has both elastic and inelastic components. The inelastic component caters to the energy requirement to operate basic household appliances like lights, fans, refrigerator etc. The elastic component takes care of more extensive needs like maintaining a pleasant temperature setting inside the household (HVAC) or charging an EV.  

\subsection{Defining the Game}
Since all stakeholders involved (consumers \& retailer) are interested in optimizing their own payoff, we formulate the problem using a game theory framework. At the start of a period, the retailer chooses a price and then consumers decide how much to consume in that period. Thus, the game $G$ is a sequential one. There is a `leader' (retailer) who takes the first turn at playing the game (by setting price) and the 'followers' (consumers) respond accordingly (by optimally deciding their consumption). Due to the hierarchy in player 'types', $G$ qualifies as a Stackelberg game which can be solved by Backward Induction. We assume that all players are rational and the game is one with complete information. Unlike the existing literature, we consider that the retailer may charge different consumers at different rates. We define the players' payoff functions in the subsequent sections.

\section{Optimal Pricing in a Grid without Distributed Generation}\label{sec:norenewable}
In this section, we assume that consumers do not have in-house renewable energy generation capabilities. So all of their demand are satisfied by the grid. First, we describe the consumer's demand model and subsequently discuss the optimal decision of a consumer given the price of a retailer (Section~\ref{sec:consume_dec}). 

To compute the optimal price, we formulate the retailer's optimization problem which turns out to be non-convex (Section~\ref{sec:retail_formulate}). So we formulate two relaxed versions of the original problem which are convex. We show that when the retailer can charge different prices to different users, she always sets higher prices to consumers who have higher willingness for consumption (Theorem~\ref{thm:fair}). Thus the proposed price mechanism induces fairness among the consumers (Section~\ref{sec:retail_fair}). We also characterize how the optimality gap can be measured between the relaxed formulations and original formulation in the non-discriminatory pricing regime (Section~\ref{sec:asymptotic_analysis}). We conclude by empirically evaluating the characteristics of the pricing scheme and its impact on the consumer's welfare and retailer revenues in  Section~\ref{sec:retail_numerical}. 

\subsection{Consumer Decision}\label{sec:consume_dec}
The net demand of a consumer has both elastic and inelastic components. We assume that the inelastic demand requirement for all consumers is the same (equal to $m^{(k)}$). The retailer charges consumers at a base rate $p_b$ for the inelastic demand $m^{(k)}$. Since the allowance is meant for satisfying basic energy requirements for a consumer, it would be unfair if the retailer chose to optimize over $p_b$. So, we assume that $p_b$ is fixed beforehand (set at a low value) and known to all. The retailer only optimizes her price for the elastic component of demand. This assumption holds throughout this section. However, \textbf{note} that our earlier assumption of all consumers having the same inelastic demand requirement is only for ease of exposition and the analysis goes through smoothly even without it.

For a given price $p_i^{(k)}$ at period $k$, consumer $i$ decides how much to consume. The decision is based on the convenience function and the price. Each consumer derives some comfort from consuming energy and we monetize it using a convenience function. The use of convenience functions is well documented in literature (\cite{samadi2010optimal}, \cite{fahrioglu1999designing} and \cite{fahrioglu2001using}). However, note that the convenience function is dependent only on elastic demand (inelastic demand does not contribute to convenience as it meets the basic requirements of a consumer). We, now, describe the convenience function and the utility of a consumer.\\

First, we introduce some notation. Let $\omega_i^{(k)}$ be the preference factor for consumer $i$ in period $k$ (measuring willingness to consume). $x_i^{(k)}$ is the elastic demand chosen by the consumer in the same period. $\alpha$ is a predetermined constant. 
\begin{assumption}\label{ass:dem_thres}
The elastic demand $x_i^{(k)}$ of consumer $i$ in period $k$ is bounded, i.e., $0 \leq x_i^{(k)} \leq \frac{\omega_i^{(k)}}{\alpha}$ $\forall$ $i$. The upper bound indicates the demand threshold beyond which convenience gets saturated and does not increase any further with increase in demand.
\end{assumption} 
The convenience function $C(.)$ must have the following desirable properties:
\begin{itemize}
    \item $C(0, \omega_i^{(k)})$ = 0, i.e., the function has a fixed point at the origin. For zero elastic demand, convenience derived must be zero.
    \item $dC(\cdot)/dx_i^{(k)} \geq$ 0. Convenience should be an increasing function of elastic demand.  
    \item $d^{2}C(\cdot)/d(x_i^{(k)})^{2} \leq$ 0. The higher the consumption of elastic demand, the lower the marginal convenience derived from it. 
    \item The convenience saturates once marginal convenience goes to zero. Thus, if a user's demand exceeds a certain threshold ($\frac{\omega_i^{(k)}}{\alpha}$), the demand will not fetch any additional convenience to the users (Assumption \ref{ass:dem_thres}). 
    \item $C(\cdot)$ is continuous and at least twice differentiable over the entire range of $x_i^{(k)}$. This is for analytical tractability. 
\end{itemize}
Taking all the above into consideration, we define our convenience function as : 
\begin{equation*}
    C(x_i^{(k)}, \omega_i^{(k)}) = \omega_i^{(k)} x_i^{(k)} - \alpha\frac{(x_i^{(k)})^{2}}{2}, \hspace{10pt}  x_i^{(k)} \leq \frac{\omega_i^{(k)}}{\alpha}
\end{equation*}
Quadratic convenience functions are quite common in the smart grid literature (\cite{samadi2010optimal}, \cite{fahrioglu1999designing}). Note that convenience function is time dependent because $\omega$ varies over time (a consumer has different willingness to consume during different times of the day). However, for simplicity, we assume that the convenience function is not correlated across time periods. We leave that as a future direction of our work.

\begin{definition}\label{defn:utility}The consumer's utility is defined as the difference between the convenience derived from the elastic demand consumption and the total price paid for the consumption. Hence, mathematically, the utility function is
\begin{equation}
    U_i^{(k)}(x_i^{(k)}, \omega_i^{(k)} | p_i^{(k)}) = C(x_i^{(k)}, \omega_i^{(k)}) - p_i^{(k)} x_i^{(k)}
\end{equation}
\end{definition}

A consumer selects $x_i^{(k)}$ which maximizes the utility. Hence, 
\begin{observation}
The optimal consumer-end elastic demand consumption in the $k^{th}$ period in response to price $p_i^{(k)}$ charged by retailer is  
\begin{align}
x_i^{k}=\max (0, \frac{\omega_i^{(k)}-p_i^{(k)}}{\alpha}).
\end{align}
\end{observation}
Consumers' demand functions which are linear in price, are very common in the smart grid literature. However, our analysis uses a piece-wise linear response function (as in Observation 1)  which is closer to reality. But it introduces non-convexity in the retailer-end optimization problem which we discuss in the following sections.

\subsection{Retailer's Decision}\label{sec:retail_formulate}
As stated earlier, the retailer charges consumers at a minimal base rate $p_b$ for the inelastic demand. $p_b$ is fixed beforehand. The retailer only optimizes price for the elastic demand. We consider a discriminatory pricing regime where different consumers are charged at different rates. 

\begin{assumption}\label{ass:omega}
The retailer knows $\omega$ (parameter indicating consumer's willingness for demand) accurately. 
\end{assumption}
\textit{Implications:} The retailer does not ask consumers to reveal their preferences ($\omega$'s). The smart meter installed in a household can analyse the user's consumption data and accurately measure $\omega$.  With the knowledge of each consumer's $\omega$, the retailer can accurately predict the optimal response of a consumer to a price. Since the consumer $i$ is not communicating their $\omega_i$, we do not run into  \textit{incentive-compatibility} issues. In scenarios where consumers would be allowed to report their $\omega_i$, they could choose to report incorrect (or lower) values deliberately to game the system and get charged at lower rates. However, we do not run into such problems here. 

\paragraph{Retailer's Objectives}: The primary objective of the retailer is to maximize their revenue. They also incur a cost associated with satisfying consumer demand. The cost is considered to be quadratic in its argument which is consistent with the existing literature (\cite{zivic2012simplified, huang2014using}).

Since electricity is an essential commodity, in most cases, because of regulatory intervention, the retailer is required to take consumer welfare into account in addition to revenue considerations. Retailers may incur penalties for setting high prices to ensure that consumers are able to consume as close as possible to their optimal levels. Note from Observation 1 that the optimal demand at price $p$ is $max(0, \frac{\omega-p}{\alpha})$. Hence, the  optimal demand is $\frac{\omega}{\alpha}$ which is realized when $p = 0$. \textit{Note that the retailer sets prices for only one period at a time and the price is announced at the start of the period.}

Thus we have the following optimization problem for the retailer:\\ \\
\textbf{Formulation 0}
\begin{eqnarray}
\text{maximize }  e_1\cdot( \sum_{i}p_i^{(k)}x_i^{(k)}) - \frac{e_2}{N} \cdot (\sum_{i} x_i^{(k)})^{2} - e_3\cdot \sum_{i} (\frac{p_i^{(k)}}{\alpha})^{2}\nonumber\\
\text{subject to }  x_i^{(k)} = \max(0,\frac{\omega_i^{(k)}-p_i^{(k)}}{\alpha})\hspace{6pt} \forall \hspace{3pt} i\nonumber\\
-\eta \leq p_i^{(k)}-p_j^{(k)} \leq \eta\hspace{6pt} \forall \hspace{3pt} i \neq j\label{eq:dis}\\
p_b \leq p_i^{(k)} \leq P \hspace{6pt} \forall \hspace{3pt} i\label{eq:lim}
\end{eqnarray}
$e_1, e_2, e_3$ are weight factors. They must be chosen judiciously depending on the need. The first term in the objective corresponds to the revenue, the second term corresponds to the cost of serving the consumption. As $N$ becomes large, the second term also becomes very large due to its quadratic nature, so we include an additional $\frac{1}{N}$ to make it comparable to the first and third terms at all times. The third term in the objective represents a penalty for setting a high price. The first constraint denotes the fact that user's consumption is given by the expression in Observation 1. The second constraint denotes that the level of discrimination is limited to $\eta$. The last constraint bounds the decision variable price. Note that $p_i^{(k)} \geq p_b$ because elastic demand should be charged at a higher rate than inelastic demand. 

\paragraph{Limiting the Discrimination}: Note that the prices charged to two different consumers differ by at most $\eta$. If $\eta=0$, we revert to the scenario where there is no discrimination. On the other hand, if we have $\eta=\infty$, we revert to the scenario where the discrimination level is not bounded at all. $\eta$ is a policy choice for the social planner. We numerically show the impact of $\eta$ on each of the retailer objectives. \\ \\
\textbf{Formulation 0} is not convex since the first constraint is a non-linear equality constraint. Thus, it is difficult to obtain an optimal price. Therefore, we propose a couple of reformulations of the original optimization problem which are convex.

\subsubsection{Formulation 1}
\begin{equation*}
        \begin{aligned}
            \underset{p_i^{(k)}, x_i^{(k)}}{\mathrm{Max}} \hspace{10pt}  e_1\cdot( \sum_{i}p_i^{(k)}x_i^{(k)}) - \frac{e_2}{N}\cdot (\sum_{i} x_i^{(k)})^{2} - e_3\cdot\sum_{i} (\frac{p_i^{(k)}}{\alpha})^{2}
        \end{aligned}
\end{equation*}
   
\begin{equation}
        \begin{split}
        \mathrm{Subject~ to:~}
        ~&x_i^{(k)} = \frac{\omega_i^{(k)}-p_i^{(k)}}{\alpha} \hspace{6pt} \forall \hspace{3pt} i\label{eq:eq}\\
        &(\ref{eq:dis})-(\ref{eq:lim})\quad x_i^{(k)} \geq 0 \hspace{6pt} \forall \hspace{3pt} i
        \end{split}
\end{equation}
In this formulation, we do away with the max term of the first constraint. This means that the equality constraint is now linear and the optimization problem is convex. Note that we have also introduced a constraint of the form $x_i^{(k)} \geq 0$ to compensate for the relaxation made earlier. This further restricts the price. By substituting $x_i^{(k)}$ in the objective, we can reformulate the same problem as a convex optimization problem in variables $\{p_i^{(k)}\}$.
\begin{equation*}
        \begin{aligned}
            \underset{p_i^{(k)}}{\mathrm{Max}} \hspace{10pt}  e_1\cdot \sum_{i}p_i^{(k)}(\frac{\omega_i^{(k)}-p_i^{(k)}}{\alpha}) - \frac{e_2}{N}\cdot (\sum_{i} \frac{\omega_i^{(k)}-p_i^{(k)}}{\alpha})^{2} - e_3 \cdot \sum_{i} (\frac{p_i^{(k)}}{\alpha})^{2}
        \end{aligned}
    \end{equation*}
   
\begin{equation}
    \begin{split}
        \mathrm{Subject\hspace{3pt} to:} 
        &(\ref{eq:dis})-(\ref{eq:lim}), \quad p_i^{(k)} \leq \omega_i^{(k)}   \forall \hspace{3pt} i
    \end{split}
\end{equation}

\subsubsection{Formulation 2}
 Instead of putting a hard constraint like $x_i^{(k)}\geq 0$, we rather put a penalty for $x_i^{(k)} < 0$. As a result, this formulation does not restrict the price like in \textit{Formulation 1}. \textit{Formulation 2} is given by : 
\begin{equation*}
       \begin{aligned}
           \underset{p_i^{(k)}, x_i^{(k)}}{\mathrm{Max}} \hspace{10pt}  e_1\cdot( \sum_{i}p_i^{(k)}x_i^{(k)}) - \frac{e_2}{N}\cdot (\sum_{i} x_i^{(k)})^{2} - e_3\cdot\sum_{i} (\frac{p_i^{(k)}}{\alpha})^2 +\\ \sum_{i} \min(0, \frac{\omega_i^{(k)} - p_i^{(k)}}{\alpha})
        \end{aligned}
\end{equation*}
   
\begin{equation}
      \begin{split}
          \mathrm{Subject\hspace{3pt} to:} 
          &(\ref{eq:dis})-(\ref{eq:lim}), x_i^{(k)} = \frac{\omega_i^{(k)}-p_i^{(k)}}{\alpha} \hspace{6pt} \forall \hspace{3pt} i
        \end{split}
\end{equation}
Note that it is convex. However, one drawback of \textit{Formulation 2} is that it does not always give the correct value of elastic load. The actual value needs to be recomputed using \textit{Observation 1} after the optimal price has been found. We can simplify this formulation further by substituting $x_i^{(k)}$ and adding a dummy variable $t_i^{(k)}$ and reformulate it as :
\begin{equation*}
        \begin{aligned}
            \underset{p_i^{(k)}, t_i^{(k)}}{\mathrm{Max}} \hspace{10pt}  e_1\cdot \sum_{i}p_i^{(k)}(\frac{\omega_i^{(k)}-p_i^{(k)}}{\alpha}) - \frac{e_2}{N}\cdot (\sum_{i} \frac{\omega_i^{(k)}-p_i^{(k)}}{\alpha})^{2} \\- e_3 \cdot \sum_{i} \left(\frac{p_i^{(k)}}{\alpha}\right)^{2} + \sum_{i} t_i^{(k)}
        \end{aligned}
\end{equation*}
   
\begin{equation}
        \begin{split}
            \mathrm{Subject\hspace{3pt} to: \hspace{3pt}} 
            &t_i^{(k)} \leq 0 \hspace{6pt} \forall \hspace{3pt} i \\
            &t_i^{(k)} \leq \frac{\omega_i^{(k)}-p_i^{(k)}}{\alpha} \hspace{6pt} \forall \hspace{3pt} i \\
            & (\ref{eq:dis})-(\ref{eq:lim})
        \end{split}
\end{equation}

\begin{theorem}\label{thm:price_comparisons}
The optimal price given by formulation 2 is always at least as high as formulation 1 for any $\eta \geq 0$. 
\end{theorem}
The result of theorem \ref{thm:price_comparisons} is intuitive. When $\eta > 0$, observe that the solution space of formulation 1 is a subset of the solution space of formulation 2. This is because of the additional constraint $x_i^{(k)} \geq 0$ in formulation 1. This implies that the optimal solution of formulation 1 is always at least a feasible solution for formulation 2. Since formulation 2 allows higher prices than formulation 1, the optimal price vector for formulation 2 can charge higher prices than formulation 1 for some consumers.  

When $\eta = 0$, the hard constraint that $x_i^{(k)} \geq 0$ $\forall$ i in formulation 1, forces the price $p^{(k)}$ to be upper-bounded by $min_{i}\{\omega_i^{(k)}\}$. On the contrary, formulation 2 only imposes a penalty whenever a user is unable to consume a positive elastic demand. So, prices for formulation 2 are expected to be higher. For a detailed proof of the same, please refer to the Appendix section of the paper.  

\subsection{Fairness}\label{sec:retail_fair}
One of the drawbacks with non-discriminatory pricing (same price for all) is that  consumers who have low willingness for demand (often stemming from poor economic conditions) pay at the same rate as that the ones who have higher willingness for demand.  Many of them are forced to manage with the inelastic demand allowance provided by the retailer as they are unable to afford any extra energy for meeting their elastic demands.  In this segment, we show that discriminatory pricing is a \textit{fair} pricing strategy. It charges higher prices to consumers with higher willingness for demand, thus lowering their consumption while consumers with low willingness for demand are charged lower rates which enables them to consume more. Hence, our proposed pricing mechanism facilitates a more equitable distribution of energy. 

\begin{theorem}\label{thm:fair}
Discriminatory pricing charges a higher price to a consumer with higher willingness for demand. Therefore, $\omega_i > \omega_j$ $\implies$ $p_i \geq p_j$. Additionally, if $\omega_i = \omega_j$, $p_i = p_j$.
\end{theorem}
\textit{Proof Sketch :} Due to space constraints, we provide a sketch of the original proof to help the reader with intuitive understanding. Note that the superscripts corresponding to period $k$ are being dropped for ease of notation.
We complete the proof by contradiction. Let $(p_1, p_2,...p_N)$ be the optimal price vector that maximizes retailer objective for a given $\eta$. We start by making a claim that $\exists$ $i$ and $j$ such that $\omega_i > \omega_j$ and yet $p_i < p_j$. However, we find that by interchanging $p_i$ and $p_j$ (i.e., charging $p_i$ to consumer $j$ and $p_j$ to consumer $i$), we are able to achieve a higher retailer objective for both formulations 1 and 2. This contradicts our claim that $p_i < p_j$ and yet the price vector is optimal. Therefore, for $\omega_i > \omega_j$, $p_i \geq p_j$. \\ \\
We follow a similar approach in order to prove the equality. Let $(p_1, p_2,...p_N)$ be the optimal price vector. We claim that $\exists$ $i$ and $j$ such that $\omega_i=\omega_j$ and yet $p_i \neq p_j$. Without loss of generality, we take $p_i > p_j$. However, we show that a higher retailer objective is achievable for both formulations 1 and 2 by charging the same price $\frac{p_i+p_j}{2}$ to both consumers $i$ and $j$, Hence, our claim cannot be true and $\omega_i = \omega_j$ $\implies$ $p_i = p_j$. \\
For a detailed proof, check the Appendix section of the paper. 

\subsection{Optimal level of discrimination $\eta$ :} In addition to the price, we can also choose to make $\eta$ a decision variable to determine the optimal level of discrimination. Both formulations 1 and 2 still continue to be convex when $\eta$ is a decision variable. We observe that the optimal $\eta = \eta^{*}$ is the critical level of discrimination at which no consumer is forced to forego their elastic demand due to high prices. At this point, the retailer revenue achieves it maximum value and revenues and elastic loads plateau beyond this point as the price vector remains unchanged. Therefore, $\eta = \eta^{*}$ signifies the limiting level of discrimination at which the energy redistribution process is complete and the maximal social benefits are visible. For formulations 1 and 2, $\eta^{*}$ is given by:
\begin{equation}
    \eta^{*} = \frac{e_1 \alpha (\omega_{L}-\omega_{S})}{2(\alpha e_1 + e_3)}
\end{equation}
where $\omega_{L} = \max_{i=1(1)N}\{\omega_i^{k}\}$ and $\omega_{S} = \min_{i=1(1)N}\{\omega_i^{k}\}$. This can be obtained by using the Karuhn-Kush-Tucker(KKT) conditions on the constrained retailer-end optimization problem.

\subsection{Prices in the Asymptotic Sense}\label{sec:asymptotic_analysis}
In Section \ref{sec:retail_formulate}, we have proposed two convex formulations for the original retailer's optimization problem. In this segment, we characterize how the prices differ among different formulations in the asymptotic regime (by taking $N\rightarrow \infty$) when $\eta=0$. We abuse the notation slightly throughout this Section \ref{sec:asymptotic_analysis} and denote the price by $p$ after dropping the superscript $(k)$ for the period. The subscript for consumer index does not appear naturally because we are in a non-discriminatory setting, hence, we also remove subscript corresponding to a user in the notation for price. 

\subsubsection{Asymptotic Analysis for Formulation 0}
In the asymptotic limit, we consider the following as objective
\begin{align}
    \lim_{N\to \infty}\frac{1}{N}f_0(p)
\end{align}
where $f_0(p)$ is the objective value for formulation $0$. In order to determine the value of the objective in the asymptotic sense, we assume that the consumer $\omega$'s are uniformly distributed in the range [$\omega_{min}$, $\omega_{max}$] throughout this section. The above assumption enables us to write by the law of large numbers
\begin{equation}\label{eq: asymp_exp}
    \lim_{N\to \infty}\frac{1}{N}f_0(p)\rightarrow \mathbbm{E}[f_0(p)]
\end{equation}
where the expectation is taken over the distribution of the $\omega$. 

Note that $f_0(\cdot)$ consists of 3 terms. Using the linearity property of expectations, finding the expectation of each of the terms would suffice.   

The analysis has to be done separately for the two intervals [$p_b$, $\omega_{min}$] and [$\omega_{min}$, $P$] since when the price is less than $\omega_{min
}$ the second term in the objective will always be positive for each user. When the price is restricted in the interval [$p_b$, $\omega_{min}$],  the optimal price $p_0^{l}$ can be found as the solution of the following constrained optimization problem :
\begin{equation*}\label{eq:f11}
    \begin{multlined}
     {\max_{p}} \hspace{3pt} g_1(p) = e_1\cdot p\frac{(\bar \omega-p)}{\alpha} - e_2 \cdot\frac{(\bar \omega-p)^{2}}{\alpha^2} - e_3 \cdot \frac{p^2}{\alpha^2}\\
    subject\hspace{3pt}to : p_b \leq p \leq \omega_{min}
    \end{multlined}
\end{equation*}
where $\bar{\omega}$ is the expected value of $\omega$.

When the price is restricted to the interval [$\omega_{min}$, $P$], the optimal price $p_0^{r}$ can be found as the solution as the following constrained optimization problem :
\begin{equation*}\label{eq:f0}
    \begin{multlined}
     {\max_{p}} \hspace{3pt} g_0(p) = e_1\cdot \frac{p(\omega_{max}-p)^{2}}{2\alpha(\omega_{max}-\omega_{min})} - e_2 \cdot  \frac{(\omega_{max}-p)^{4}}{4\alpha^{2}(\omega_{max}-\omega_{min})^{2}} - e_3 \cdot \frac{p^2}{\alpha^2}\\
    subject\hspace{3pt}to :\omega_{min} \leq p \leq P
    \end{multlined}
\end{equation*}
It is important to note that $g_0(\omega_{min}) = g_1(\omega_{min})$ which implies that formulation 0 is continuous in the entire interval [$p_b$, $P$] in the asymptotic sense. This leads us to the following observation : 
\begin{observation}\label{obs:form0}
The optimal price $p_0^{*}$ proposed by formulation 0 in the asymptotic sense can be obtained as follows :
\begin{equation}
    \begin{split}
        p_0^{*} &= p_0^{l}, \hspace{10pt}if \hspace{3pt} g_1(p_0^{l}) > g_0(p_0^{r})\\
                &= p_0^{r}, \hspace{10pt}if \hspace{3pt} g_1(p_0^{l}) \leq g_0(p_0^{r})\\
    \end{split}    
\end{equation}
\end{observation}

\subsubsection{Asymptotic Analysis for Formulation 1}
In formulation 1, we use $x_i = \frac{\omega_i-p}{\alpha}\hspace{3pt} \forall\hspace{3pt} i$. Therefore, $\mathbbm{E}(x_i) = \frac{\bar \omega - p}{\alpha}$ and $\text{Var}(x_i) = \frac{\text{Var}(\omega)}{\alpha^{2}}$ which is independent of $p$. Thus, in accordance with Eq.\ref{eq: asymp_exp}, we have the following observation:

\begin{observation}\label{obs:form1}
The optimal price $p_1^{*}$ proposed by formulation 1 in the asymptotic sense, can be obtained by solving the following constrained optimization problem in $p$:
\begin{equation*}
    \begin{multlined}
     {\max_{p}} \hspace{3pt} e_1\cdot p\frac{(\bar \omega-p)}{\alpha} - e_2 \cdot\frac{(\bar \omega-p)^{2}}{\alpha^2} - e_3 \cdot \frac{p^2}{\alpha^2}\\
    subject\hspace{3pt}to : p_b \leq p \leq \omega_{min}
    \end{multlined}
\end{equation*}
\end{observation}
The original form of the constraint should have been $p_b \leq p \leq \min_{i=1,\ldots,N} \{\omega_i\}$. However, as $N \to \infty$, we must consider the optimization problem where the constraint is satisfied for every instance of $\omega$, hence, the upper bound becomes $\omega_{min}$. Thus, the constraint takes the form $p_b \leq p \leq \omega_{min}$.\\ \\
\textit{Uniqueness : }Note that the objective function is concave in $p$ and so it has a unique maxima. However, the maxima point may not be feasible. Let the maxima obtained from the first order condition be denoted by $p^{*}$ for the objective considering the unconstrained problem. $p^{*}$ is given by $p^{*} = \frac{\bar \omega(\alpha e_1 + 2e_2)}{2(\alpha e_1 + e_2 + e_3)}$. Now, the optimum of the constrained problem, i.e., of the problem described in Observation~\ref{obs:form1} is given by :
\begin{equation}
   \begin{split}
        p_1^{*} & = \omega_{min}, \hspace{10pt}if\hspace{3pt}p^{*} > \omega_{min}\\
                & = p^{*}, \hspace{10pt}if\hspace{3pt}p_b < p^{*} \leq \omega_{min}\\
                & = p_b, \hspace{10pt}if\hspace{3pt}p^{*} \leq p_b
   \end{split} 
\end{equation}

\subsubsection{Asymptotic Analysis for Formulation 2}
If we recall formulation 2, the only problematic term in the objective is the load error term. Therefore, we need to find the expected value of the error in load calculation. The analysis has to be done separately for the two intervals [$p_b$, $\omega_{min}$] and [$\omega_{min}$, $P$]. For the first interval [$p_b$, $\omega_{min}$], formulation 2 reduces to formulation 1 and the optimal price $p_2^{l}$ can be found as the solution of the following constrained optimization problem :
\begin{equation*}\label{eq:f1}
    \begin{multlined}
     {\max_{p}} \hspace{3pt} g_1(p) = e_1\cdot p\frac{(\bar \omega-p)}{\alpha} - e_2 \cdot\frac{(\bar \omega-p)^{2}}{\alpha^2} - e_3 \cdot \frac{p^2}{\alpha^2}\\
    subject\hspace{3pt}to : p_b \leq p \leq \omega_{min}
    \end{multlined}
\end{equation*}
For the second interval [$\omega_{min}$, $P$], the optimal price $p_2^{r}$ can be found as the solution as the following constrained optimization problem :
\begin{equation*}\label{eq:f2}
    \begin{multlined}
     {\max_{p}} \hspace{3pt}g_2(p) = e_1\cdot p\frac{(\bar \omega-p)}{\alpha} - e_2 \cdot \frac{(\bar \omega-p)^{2}}{\alpha^2}  - e_3 \cdot \frac{p^2}{\alpha^2} -\frac{(p-\omega_{min})^{2}}{2\alpha(\omega_{max}-\omega_{min})}\\
    subject\hspace{3pt}to : \omega_{min} \leq p \leq P
    \end{multlined}
\end{equation*}
It is important to note that $g_1(\omega_{min}) = g_2(\omega_{min})$ which implies that formulation 2 is continuous in the entire interval [$p_b$, $P$] in the asymptotic sense. This leads us to the following observation : 
\begin{observation}\label{obs:form2}
The optimal price $p_2^{*}$ proposed by formulation 2 in the asymptotic sense can be obtained as follows :
\begin{equation}
    \begin{split}
        p_2^{*} &= p_2^{l}, \hspace{10pt}if \hspace{3pt} g_1(p_2^{l}) > g_2(p_2^{r})\\
                &= p_2^{r}, \hspace{10pt}if \hspace{3pt} g_1(p_2^{l}) \leq g_2(p_2^{r})\\
    \end{split}    
\end{equation}
\end{observation}
\textit{Uniqueness : } It can be shown that the optimal price for formulation 2 given by $p_2^{*}$ is \textit{unique}. For the detailed proof, please refer to the Appendix (Section \ref{sec:appendix}).

From the asymptotic analysis, it is observed that formulation 2 proposes a price at least as high as formulation 1. Actually, the result holds even in the non-asymptotic case (when $N$ is finite), which we have already shown in Theorem \ref{thm:price_comparisons}. 

Now that we have analytically found the prices proposed by each of the formulations in the asymptotic sense, we can also numerically evaluate the optimality gaps across them. The analysis can be found in Section \ref{sec:retail_numerical}.

\subsection{Numerical Experiments}\label{sec:retail_numerical}
In this segment, we provide numerical insights for the formulations that have been developed above.
\subsubsection{Simulation Setup}
We assume $\alpha = 2$ across all consumers and across all periods of the day. There are 6 periods, each of duration 4 hours with the first period starting at midnight. The inelastic demands ($m^{(k)}$ values) are assumed to be identical for all consumers. Inelastic demand values across 6 periods are as follows : [0.16, 0.39, 0.63, 0.51, 0.78, 0.52] units which also represents that during the peak time the base demand is also high. Consumer $i$'s preference parameter  $\omega_i^{(k)}$  for elastic demand for period $k$ is assumed to satisfy the following relationship $\omega_i^{(k)} = (0.75 + 0.5m_k)\times \omega_i \hspace{3pt}\forall\hspace{3pt}i$ where $\omega_i\sim \mathcal{U}[3,7]$ drawn independently from other consumers. During periods of peak load, $m^{(k)}$ will be higher which means that consumer willingness for demand will also be higher. For all our experiments, we assume a base price $p_b = Re.1/unit$ unless otherwise specified.

\subsubsection{Non-discriminatory Setting}
$\eta = 0$ translates to the case when the pricing is non-discriminatory in nature. The same price is charged to all users in a given period. In this segment, we will investigate the effect of weights $e_1$, $e_2$ and $e_3$ on the final prices, retailer revenues, elastic load and consumer welfare. We will also compare results across our relaxed formulations to identify the variation of prices. {We fix $e_3 = 1$ and choose $e_1, e_2$ such that $e_1 + \frac{e_2}{\beta} = 1$. $\beta$ is a pre-determined constant that makes $\frac{e_2}{\beta}$ lie between $0$ and $1$. This allows us to observe the nature of the variation as the function of a single variable $e_1$. We have set $\beta = 10$ throughout the analysis.}\\
    \begin{figure}[h]
    \centering
    \includegraphics[width = 0.9\textwidth]{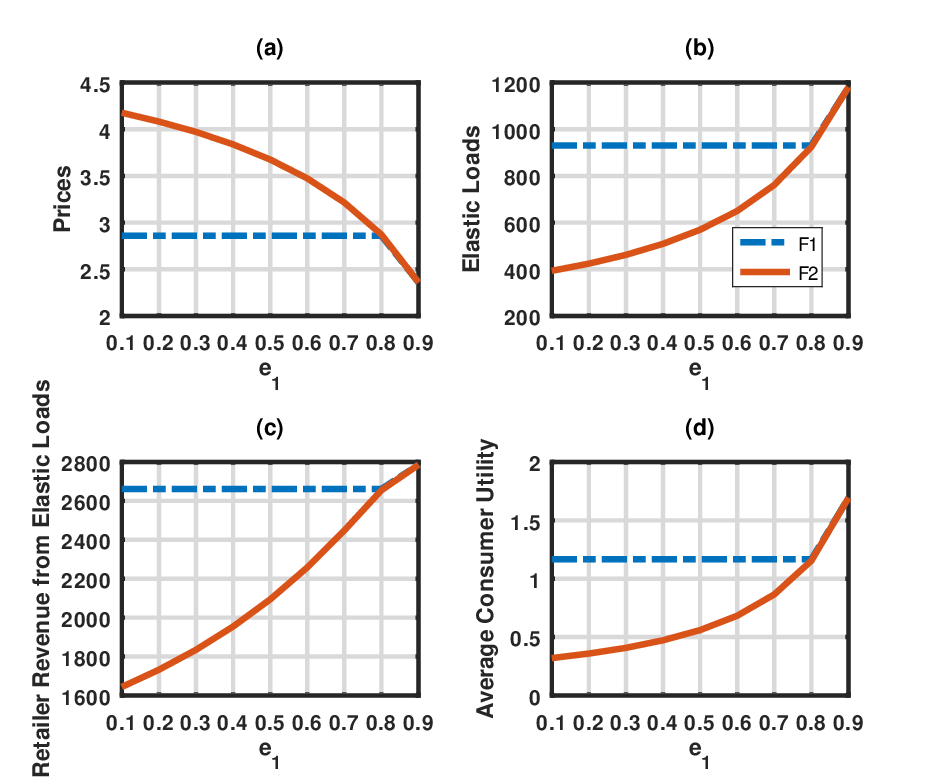}
    \vspace{-0.2in}
    \caption{Variation of key metrics across formulations F1 \& F2 as a function of $e_1$}
    \label{fig:e1_variation}
    \vspace{-0.15in}
    \end{figure} 
As we increase $e_1$, prices decrease but elastic demand consumption increases (Fig.~\ref{fig:e1_variation}). This leads to overall increase in retailer's revenues and average consumers' utilities. The reason is that when $e_1$ increases, the retailer puts a higher priority on maximizing revenue. Also an increase in $e_1$ leads to a decrease in $e_2$ (priority on reducing load) since we have forced $e_1 + \frac{e_2}{\beta} = 1$. Hence, the retailer would set lower prices to increase consumption and increase the revenue. The impacts of $e_2$ and $e_3$ are straightforward. Higher $e_2$ would lead to higher prices in order to decrease the total load. Similarly for higher $e_3$, the retailer would try to maximize the consumer's utility and reduce the price. Formulation 1 gives constant prices and constant loads/revenues for low values of $e_1$ because of the hard-constraint that $x_i^{(k)} \geq 0$ $\forall$ $i$. Hence, the price gets constricted at $p = min_i\{\omega_i^{(k)}\}$. At high values of $e_1$ (which means low values of $e_2$), the load error term in formulation 2 vanishes and both formulations give the same price. Hence, in such cases, it is convenient to use formulation 1. Also observe that formulation 2 always predicts a price which is at least as high as formulation 1.\\ 
    \begin{figure}[h]
    \centering
    \includegraphics[width = 0.8\textwidth]{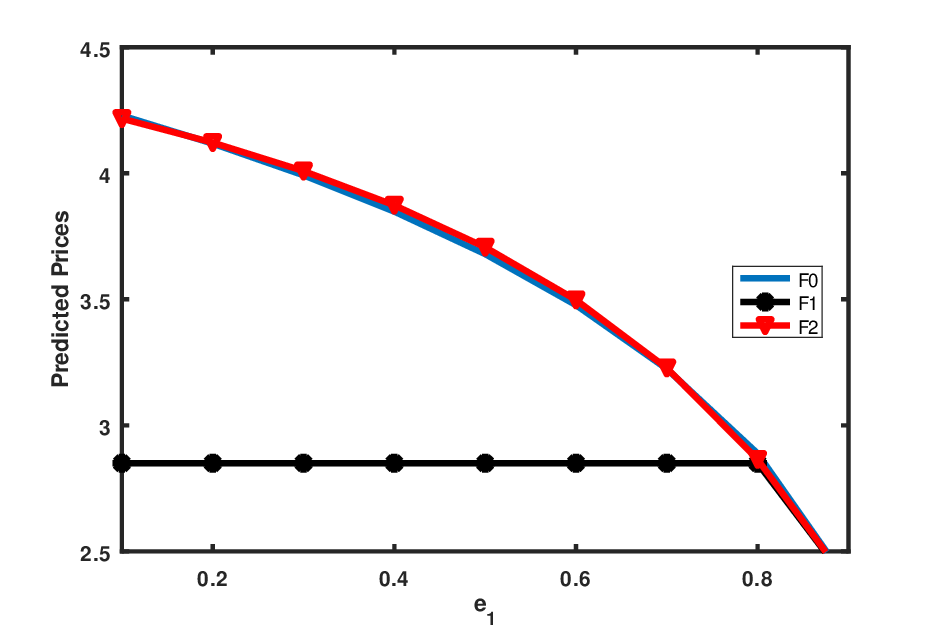}
    \vspace{-0.2in}
    \caption{Optimality gap between F0 and relaxed formulations F1 and F2 in asymptotic sense}
    \label{fig:asymptotic}
    \vspace{-0.15in}
    \end{figure} 

When number of consumers in the community is assumed to be infinitely large (asymptotic sense), we show how prices computed by relaxed formulations 1 and 2 vary with weights. From Formulation 1, we obtain  a price equal to $\omega_{min}$ for most cases, while from formulation 2, the proposed price closely tracks the price proposed by the original formulation, i.e., formulation 0. One immediate inference from Fig. \ref{fig:asymptotic} is that formulation 2 is the better of the two relaxations in terms of how closely it can approximate the original. However, at higher values of $e_1$, all the three formulations converge. 

\subsubsection{Discriminatory Setting (Positive $\eta$)}
When $\eta > 0$, our pricing model becomes discriminatory in nature, charging different prices to different users. We make several interesting observations as we vary $\eta$. \\
    \begin{figure}[h]
    \centering
    \includegraphics[width = 0.9\textwidth]{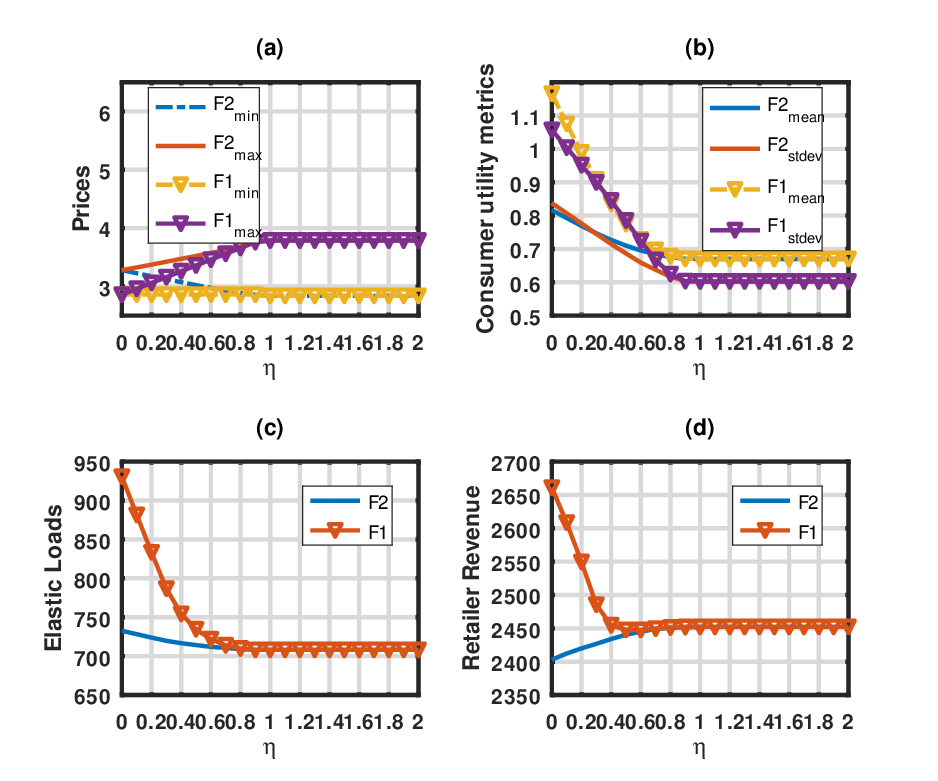}
    \vspace{-0.2in}
    \caption{Variation of key metrics for formulations 1 \& 2 as a function of $\eta$}
    \label{fig:eta_variation}
    \vspace{-0.15in}
    \end{figure} 
    \begin{figure}[h]
    \centering
    \includegraphics[width = 0.9\textwidth]{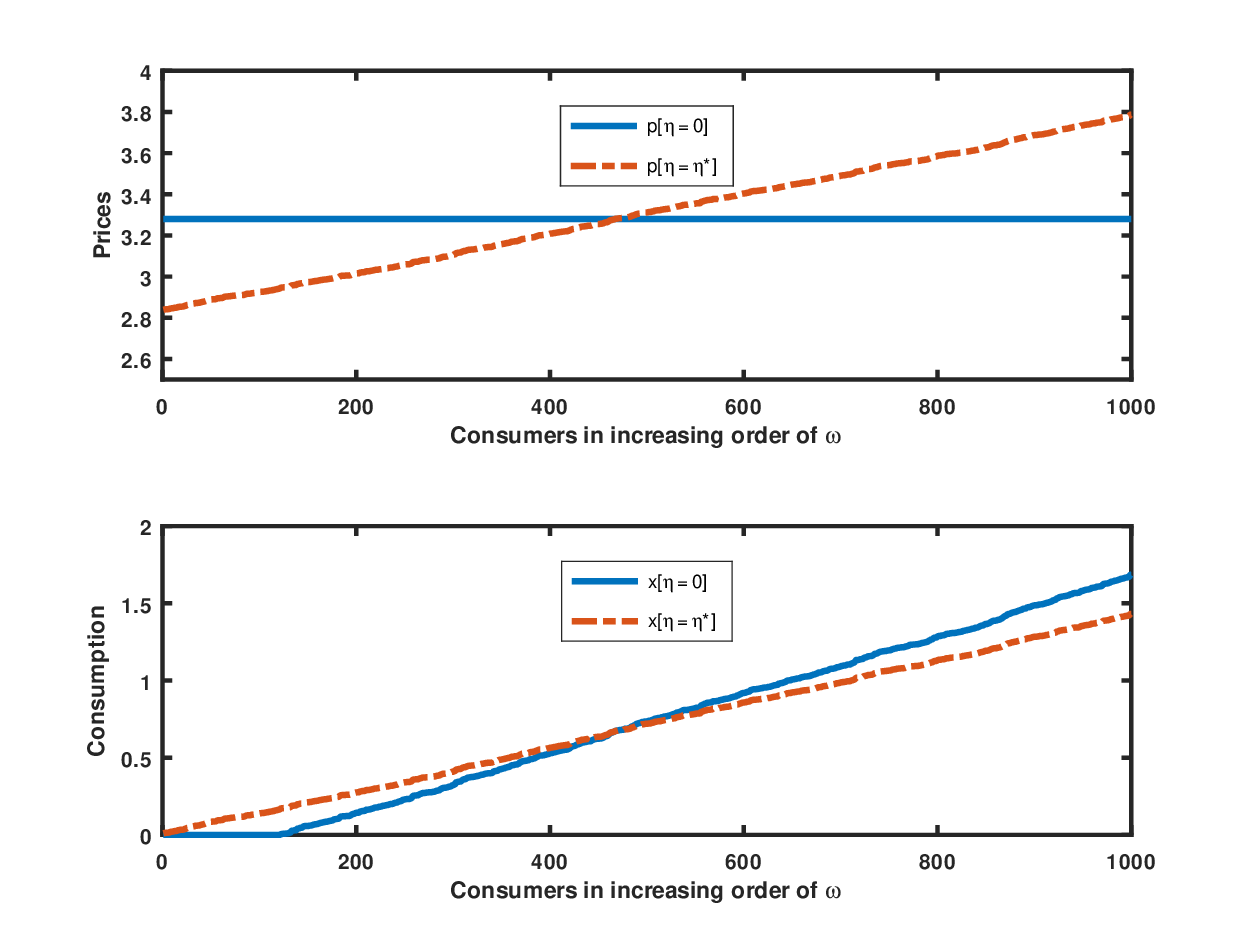}
    \vspace{-0.2in}
    \caption{Energy redistribution effects of discriminatory pricing}
    \label{fig:redistribution}
    \vspace{-0.15in}
    \end{figure} 
Discrimination affects the two relaxations in different ways. Formulation 1 proposes a small price in the non-discriminatory scenario in order to ensure that all consumers have positive elastic demands. This leads to very high elastic loads for the retailer. As the level of discrimination increases, the elastic load is reduced  which also decreases the retailer's revenue. However, the retailer's objective increases as the cost for fulfilling the demand decreases. 

For formulation 2, with an increase in $\eta$, retailer's revenue increase gradually (Fig.~\ref{fig:eta_variation}), which means that this pricing scheme is lucrative to her. Elastic load remains roughly constant even when $\eta$ varies (Fig.~\ref{fig:eta_variation}). Average consumer's utility decreases with an increase in $\eta$ (Fig.~\ref{fig:eta_variation}). Hence, it may appear that increase in $\eta$ hurts the total consumer's utility although it favours the retailer. However, it is important to note that the standard deviation in elastic demand consumption across the consumer base decreases as we increase $\eta$. The users' who have lower willingness for demand, they consumer more (Fig.~\ref{fig:redistribution}) and who have higher willingness for demand consume less as $\eta$ increases which is also suggested by Theorem~\ref{thm:fair}. Thus, {\em discriminatory pricing leads to a fairer distribution of energy in the community where high-end consumers no longer have an upper hand}. Generally, high consuming users are often related with higher wealth, thus, discriminatory price can help in eliminating the social inequalities. 

\begin{figure}[!ht]
    \centering
    \includegraphics[width = 0.8\textwidth]{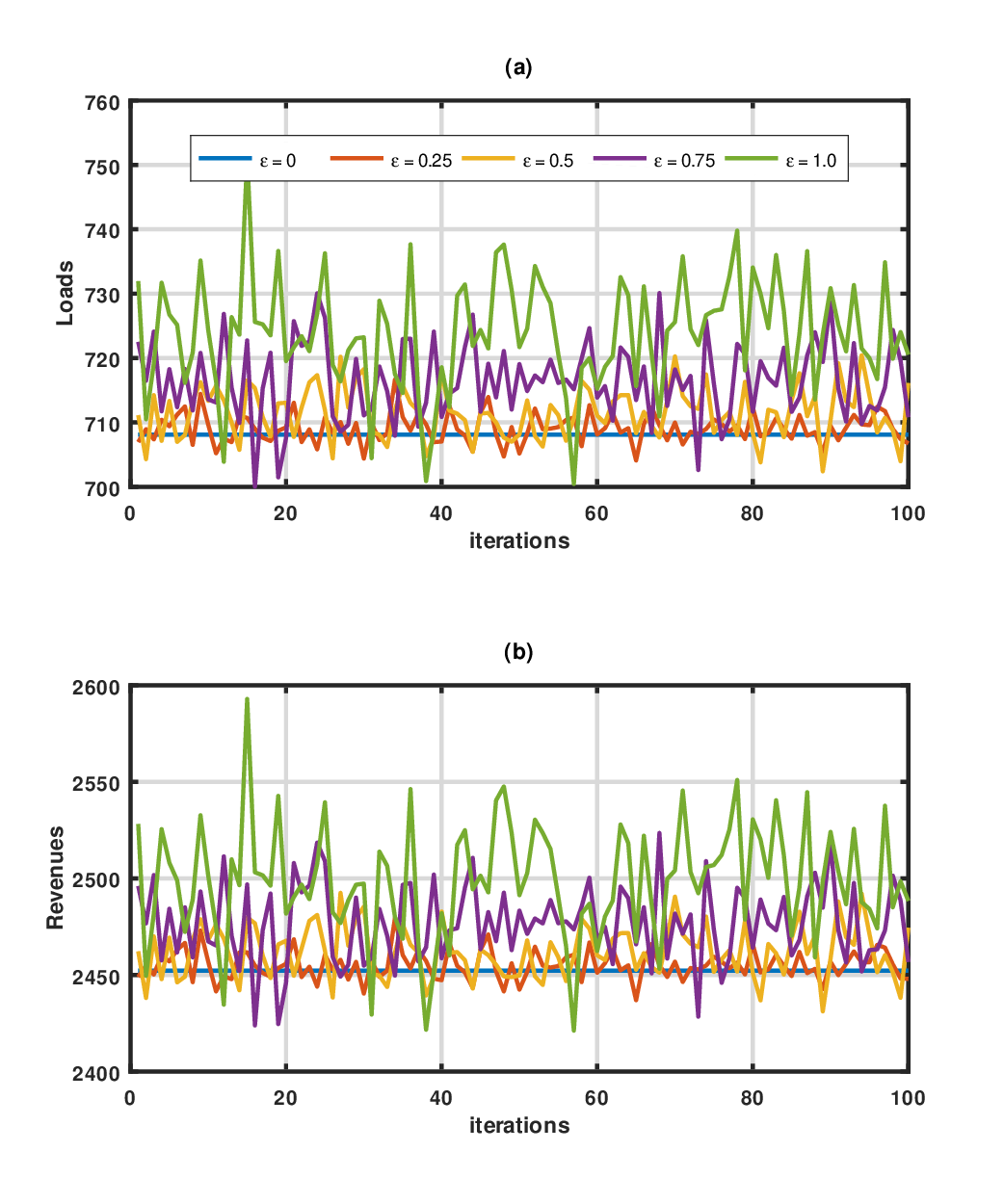}
    \vspace{-0.2in}
    \caption{Sensitivity Analysis : All results are at $\eta = \eta^{*}$, so both formulations 1 and 2 propose identical prices.}
    \label{fig:sensitivity}
    \vspace{-0.15in}
\end{figure}

\subsubsection{Sensitivity Analysis}
We have assumed earlier (Assumption \ref{ass:omega}) that the retailer has complete information about each consumer's willingness for demand parameter $\omega_i^{(k)}$. However, that retailer may not have the exact information, rather, it may have an estimate of $\omega_i^{(k)}$.  In the following, we empirically evaluate how the error in estimation impacts the retailer's revenues. \\

Let $\bar \omega_i^{(k)}$ represents the estimated (by retailer) value of $\omega_i^{(k)}$ (elastic demand willingness parameter for consumer $i$ in period $k$). Let $\epsilon$ represent the level of error in  estimation of $\omega_i^{(k)}$ in real-time. Therefore, the true-value of $\omega_i^{(k)} = \bar \omega_i^{(k)} \pm \delta$ where $\delta$ $\mathtt{\sim}$ $\mathcal{U}[0, \epsilon]$. We vary $\epsilon$ from $0$ to $1$ in steps of $0.25$ and for each value of $\epsilon$, we obtain 100 simulation runs. $\eta$ is taken as $\eta^{*}$. It is observed that as $\epsilon$ increases, the average estimation error in load and revenue increases which is expected. Also, in most of the cases, the retailer is found to underestimate the elastic load and revenue values while setting prices.

\section{Optimal Pricing when Users have Renewable Energy Resources}\label{sec:renewable}
In this section, we consider the scenario where each consumer has renewable energy generation capabilities. The renewable energies can range from solar, biomass to wind energies. Note that when a user is equipped with renewable energies, it may feed back energy to the grid. We assume the popular net-metering mechanism. Thus, the energy which is fed back is compensated at the same buying price and effectively, the consumer only pays for the net energy purchased from the grid. Since a user can technically produce energy, we denote it as a prosumer (producer+consumer). 

In Section~\ref{sec:prosumer}, we describe the prosumer's utility model. In Section~\ref{sec:retailer_renew}, we describe the retailer's optimization problem. Section~\ref{sec:fair_renew} describes how discriminatory pricing induces a fairness in this regime as well. Finally, in Section~\ref{sec:numerical_renew}, we empirically evaluate our price mechanism and its impact on the prosumers' utilities and the retailer's objectives.  

\subsection{Prosumer Decision}\label{sec:prosumer}
In the $k^{th}$ period, consumer $i$ has an solar energy generation amounting to $s_i^{(k)}$. This is complemented by a purchase of amount $x_i^{(k)}$ from the retailer at rate $p_i^{(k)}$. In case, $s_i^{(k)}$ is beyond what is required in the household, it sells back $y_i^{(k)}$ at same retail rate. $z_i^{(k)}$ is the net energy transaction made, i.e., $z_i^{(k)} = x_i^{(k)} - y_i^{(k)}$. 

Recall from Definition~\ref{defn:utility} that the prosumer's utility is defined as the difference between the convenience derived from the elastic demand consumption and the price paid for the net purchase from the grid. $z_i^{(k)}+s_i^{(k)}$ is the total demand consumption by the prosumer in the $k^{th}$ period, hence $z_i^{(k)}+s_i^{(k)}-m^{(k)}$ is the corresponding elastic demand consumption. Mathematically, the  utility function is given by
\begin{equation}
    \begin{aligned}
    U_i^{(k)}(z_i^{(k)}+s_i^{(k)}-m^{(k)}, \omega_i^{(k)} | p_i^{(k)}) = C(z_i^{(k)}+s_i^{(k)}-m^{(k)}), \omega_i^{(k)}) - p_i^{(k)} z_i^{(k)}\nonumber
    \end{aligned}
\end{equation}
Recall that $m^{(k)}$ is the inelastic demand that needs to be met at any cost. Hence, similar to Observation 1, we obtain that 
\begin{observation} The net optimal grid purchase in the $k^{th}$ period in response to price $p_i^{(k)}$ set by the retailer for both retail and sell-back is given by 
\begin{align}
    z_i^{(k)} = \max \{m^{(k)}-s_i^{(k)}, m^{(k)}-s_i^{(k)}+\frac{(\omega_i^{(k)}-p_i^{(k)})}{\alpha}\}
\end{align} 
\end{observation}
Positive $z_i^{(k)}$ indicates that renewable energy generation is insufficient and purchase has to be made from the grid to meet residual demand. While negative $z_i^{(k)}$ indicates that the renewable energy generated exceeds the requirement or it is more profitable to sell-back energy by consuming less. Note that when the grid is congested, the retailer can choose higher prices to incentivize the prosumers to sell back more. 

\subsection{Retailer Decision}\label{sec:retailer_renew}
The retailer sets the price $p_i^{(k)}$ for the $i^{th}$ prosumer in the $k^{th}$ period. The same price is applicable for both purchase and sell-back. The prosumer employs pricing regime where it can select different prices to different prosumers. In the next section, we show that such a pricing regime is fair since it induces higher prices for the ones who has higher willingness for consumption from the grid. We numerically evaluate the impact of this price mechanism on the revenue of the retailer and the users' utilities.  

In addition to the assumptions in Section \ref{sec:norenewable}, we have the following:
\begin{assumption}\label{ass:renew}
The retailer can accurately predict $s_i^{(k)}$ for each prosumer $i$ at the start of each period $k$.
\end{assumption}
Note that in general, whenever a prosumer installs a solar panel, it needs to inform the utility company. Hence, the assumption that the utility company can monitor the generation of solar energy at the prosumer's premises, is a valid one. Also, in the event that prosumers were asked to report $s_i^{(k)}$, they would have an incentive to misreport (Refer Theorem \ref{thm:renewable}). Assumption \ref{ass:renew} ensures that there are no issues with respect to incentive-compatibility.

With the knowledge of $s_i^{(k)}$ and $\omega_i^{(k)}$ and the form of the convenience function described in Section~\ref{sec:prosumer}, the retailer  knows the net optimal purchase amount for that user using Observation 2. \\ \\
\textit{Retailer's Objectives :} As mentioned in Section~\ref{sec:norenewable}, the retailer will try to maximize her own revenue, minimize the cost, and maximize the user's welfare. Thus, the retailer's optimization problem is \\ \\
\textbf{Formulation 3:}
\begin{eqnarray}
\text{maximize}  e_1(\sum_{i}p_i^{(k)} z_i^{(k)}) - \frac{e_2}{N}(\sum_{i} z_i^{(k)})^{2}\nonumber - e_3\sum_{i}(\frac{p_i^{(k)}}{\alpha})^{2}\nonumber\\
\text{subject to }  z_i^{(k)} = \max(m^{(k)}-s_i^{(k)}, m^{(k)}-s_i^{(k)}+\frac{\omega_i^{(k)}-p_i^{(k)}}{\alpha})\nonumber\\
0 \leq p_i^{(k)} \leq P,\quad 0 \leq \sum_{i} z_i^{(k)}\nonumber
\end{eqnarray}
The first term in the objective corresponds to the revenue, the second term corresponds to the cost of serving the consumption. Note that even when $z_i^{(k)}$ is negative, the retailer needs to dispatch the additional energy which incurs a cost. This is because the balance needs to maintained between the supply and demand, and even when supply exceeds the demand the retailer pays a penalty for the imbalance. The third term in the objective represents a penalty for setting too high prices.  

The first constraint denotes the fact that user's consumption is given by the expression in Observation 2. The second constraint indicates that the retailer should be able to sell a net positive amount of energy to the users which will result in her revenue. The last constraint provides an upper and lower limit for the decision variable $p_i^{(k)}$. 

Formulation 3 is non-convex because the first constraint is a non-linear equality constraint.  So, we relax the constraint and reformulate the problem as a convex one. The reformulations are exactly identical in structure to the ones provided in Section \ref{sec:norenewable} and so we omit them here.

\subsection{Theoretical Results: Significance of Discriminatory pricing regime}\label{sec:fair_renew}
Let us start with a very important observation which will lead us to Thm \ref{thm:renewable}.
\begin{observation}
For the same price $p^{(k)}$, if a prosumer generates more energy, then he will sell-back more if he chooses to sell back at all.
\end{observation}
\begin{proof}
Recall that $z_i^{(k)} = Max(m^{(k)} - s_i^{(k)} + \frac{\omega_i^{(k)}-p^{(k)}}{\alpha}, m^{(k)}-s_i^{(k)})$. Negative $z_i^{(k)}$ implies sell-back by prosumer $i$. Now, there are 2 cases, case 1 where $\omega_i^{(k)} > p^{(k)}$ and case 2 where $\omega_i^{(k)} \leq p^{(k)}$. In case 1, we can trivially note that $z_i^{(k)} = m^{(k)} - s_i^{(k)} + \frac{\omega_i^{(k)}-p^{(k)}}{\alpha}$. Observe that if $z_i^{(k)} < 0$, higher $s_i^{(k)}$ indicates more negative $z_i^{(k)}$ which implies higher sell-back. Similarly, in case 2, $z_i^{(k)} = m^{(k)}-s_i^{(k)}$. If $m^{(k)} < s_i^{(k)}$ and the prosumer does choose to sell-back, higher $s_i^{(k)}$ always translates into higher sell-back. Also observe that, higher $s_i^{(k)}$ indicates lower grid purchases if $z_i^{(k)} > 0$.
\end{proof}
The above result is intuitive. When a prosumer generates more energy, he/she will try to sell-back more. 

\begin{theorem}\label{thm:renewable}
If two prosumers have the same willingness for demand, then discriminatory pricing chooses a lower price for one who produces more energy. Therefore, if $\omega_i^{(k)} = \omega_j^{(k)}$, then $s_i^{(k)} \geq s_j^{(k)}$ $\implies$ $p_i^{(k)} \leq p_j^{(k)}$ where equality between prices hold only when $s_i^{(k)} = s_j^{(k)}$.
\end{theorem}
The intuition behind theorem \ref{thm:renewable} is that discriminatory pricing benefits prosumers who contribute energy to the grid when they have same level of willingness for demand. Prosumers who have to purchase energy from the grid to meet residual demand, the greater the contribution of the prosumers, the greater discount they get on the retail price. On the other hand, when prosumers do not have to purchase from the grid, it lowers the sell back price (since the price the prosumer gets is smaller) for them to dissuade them from selling back very high amounts. Thus it prevents the issue of retailer being left over with excess sell-back energy at the end of the period.\\ \\
\textit{Proof Sketch :} For simplicity of notation, we drop the superscript $(k)$ indicating period of the day as all quantities are considered for the same period. Let us consider 2 consumers $i$ and $j$ such that $\omega_i = \omega_j$ and $s_i > s_j$. Let the optimal price vector be given by $p^{*} = (p_1, p_2,....p_n)$ where $p_i$ and $p_j$ are the prices charged to consumers $i$ and $j$ respectively. Let if possible, $p_i > p_j$.\\
Now, let us consider a slightly modified price vector $p^{'}$ such that consumer $i$ is charged the price $p_j$ and consumer $j$ is charged the price $p_i$, i.e., the prices $p_i$ and $p_j$ are interchanged. \\
We start by considering the objective in formulation 1 ($f_1$). Since $p^{*}$ is assumed optimal for formulation 1, we must have $f_1(p^{*}) \geq f_1(p^{'})$. Since $p^{*}$ and $p^{'}$ differ only in 2 components, we ignore the terms of the objective which are not affected by the switch. Therefore,
\begin{equation*}
    \begin{multlined}
    p_i(m_k-s_i+\frac{\omega_i-p_i}{\alpha}) + p_j(m_k-s_j+\frac{\omega_j-p_j}{\alpha}) \\\geq p_i(m_k-s_i+\frac{\omega_i-p_j}{\alpha}) + p_j(m_k-s_j+\frac{\omega_j-p_i}{\alpha})\\
    \implies -p_i s_i - p_j s_j \geq -p_j s_i - p_i s_j\\
    \implies p_i s_i + p_j s_j \leq p_j s_i + p_i s_j\\
    \implies p_i (s_i - s_j) \leq p_j (s_i - s_j)\\
    \implies p_i \leq p_j \hspace{10pt}(since\hspace{3pt} s_i \geq s_j)
    \end{multlined}
\end{equation*}
Thus, we arrive at a contradiction. Therefore, when $\omega_i = \omega_j$, for $s_i > s_j$, $p_i \leq p_j$. Let's try to extend the proof to formulation 2. For extension, we just need to show additionally that
\begin{equation*}
    \begin{multlined}
    m_k - s_i + Min(0, \frac{\omega_i-p_j}{\alpha}) +  m_k - s_j + Min(0, \frac{\omega_j-p_i}{\alpha}) \geq \\
    m_k - s_i + Min(0, \frac{\omega_i-p_i}{\alpha}) +  m_k - s_j + Min(0, \frac{\omega_j-p_j}{\alpha})
    \end{multlined}
\end{equation*}
Since $\omega_i = \omega_j$, equality holds (can be trivially verified). Thus, the claim holds for formulation 2 as well. The case for the equality is also very trivial and for the sake of brevity, is omitted here. 

\subsection{Numerical Experiments}\label{sec:numerical_renew}
In this segment, we provide numerical insights for the scenario where some of the consumers are prosumers and they feed energy back to the grid. A net-metering framework is assumed where sell-back is compensated at the retail rate. 
\subsubsection{Simulation Setup}
The simulation setup is exactly similar to the one we had in the last section. Additionally, we specify the specifications for the renewable energy production inhouse. We assume that household $i$ produces energy $s_i^{(k)}$ in the $k^{th}$ period of the day where $s_i^{(k)} \sim \mathcal{U}[0, a+\frac{b\omega_i^{(k)}}{\alpha}]$. We set $a = 1$ and $b = 1$. This helps us to explore the results for a community which has users from across the spectrum (from consumers who produce energy in negligible amounts to prosumers who produce much in excess of their requirements).

\subsubsection{Discriminatory Setting}
We vary the level of price-discrimination, $\eta$, and observe the effects of the same on prices, loads and revenues. One important observation is that with net-metering, both relaxations of original formulation propose identical price vectors. This is because prices are significantly lower than the previous scenario (where there is no distributed generation), causing the condition $p_i^{(k)} \leq \omega_i^{(k)}$ to be generally valid. As a result, both relaxations reduce to the same optimization problem and produce the same results.  

As $\eta$ increases, the amount of sell-back decreases and the net load for the retailer remains more or less the same. Sell-back decreases because prosumers with high sell-back are compensated at lower rates which dissuades them from selling exorbitant amounts. However, revenue increases sharply which is favourable for the retailer. This can be attributed to the fact that lower sell-back leads to lower payouts to prosumers. As grid operation costs remain unchanged due to a stable net load, revenues are bound to increase. Discrimination can be a very handy tool for the retailer in order to keep sell-back amounts in check. This is important because too much sell-back is detrimental to her revenue while too less can cause the net demands to spiral out of control. 
\begin{figure}[!ht]
    \centering
    \includegraphics[width = 0.8\textwidth]{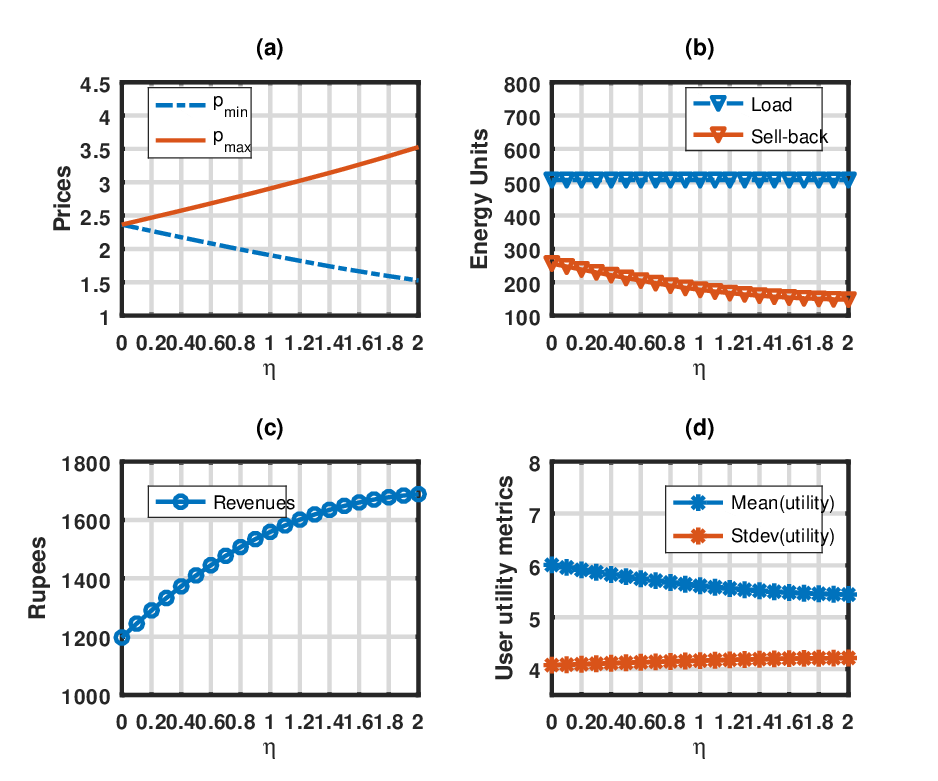}
    \vspace{-0.2in}
    \caption{Variation of key metrics as a function of $\eta$}
    \label{fig:disc_nm}
    \vspace{-0.15in}
\end{figure}

\subsubsection{Effect of varying degrees of penetration of renewable energy}
\begin{figure}[!ht]
    \centering
    \includegraphics[width = 0.8\textwidth]{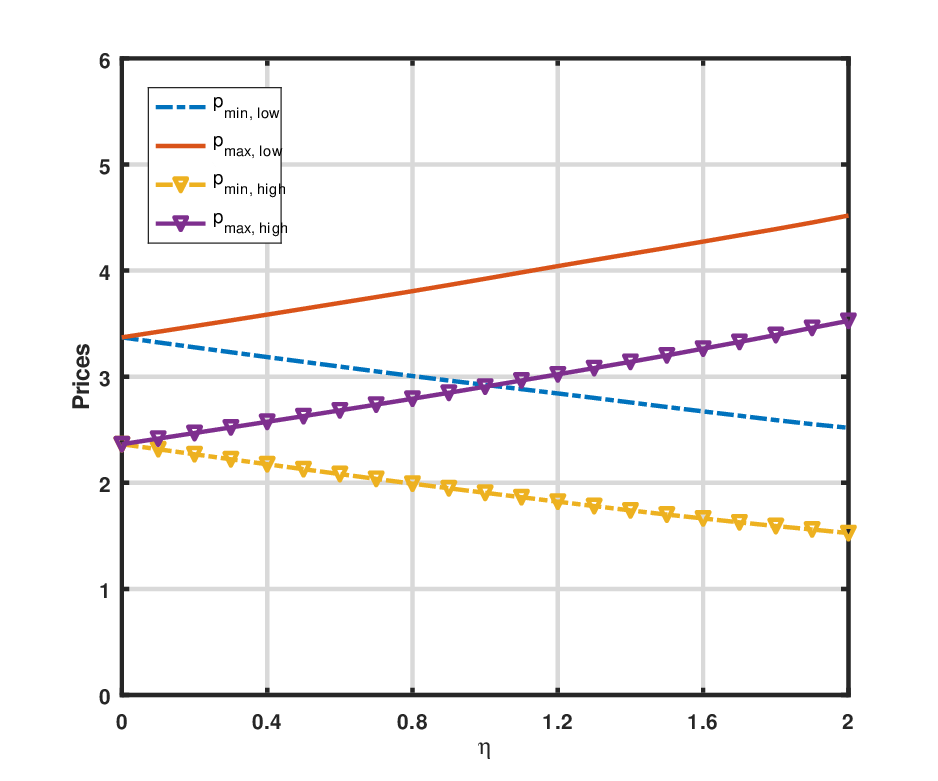}
    \vspace{-0.2in}
    \caption{Variation of prices as a function of $\eta$ with varying degrees of renewable energy penetration}
    \label{fig:penetration_nm}
    \vspace{-0.15in}
\end{figure}

Penetration of renewable energy refers to the installed capacity of renewable energy generation in the community. In this segment, we observe the effect of degree of penetration on prices. In the first set of experiments, we numerically evaluate the effect on prices under different levels of discrimination when the installed capacity is halved across all households compared to previous values (Fig \ref{fig:penetration_nm}). We observe that as the penetration decreases, the prices increase. This is intuitive because a lower penetration of renewable energy implies that the retailer's net load is significantly higher and hence the prices are higher. The results presented here are for two cases $high$ and $low$, where $low$ corresponds to $50\%$ penetration. \\

\begin{figure}[!ht]
    \centering
    \includegraphics[width = 0.8\textwidth]{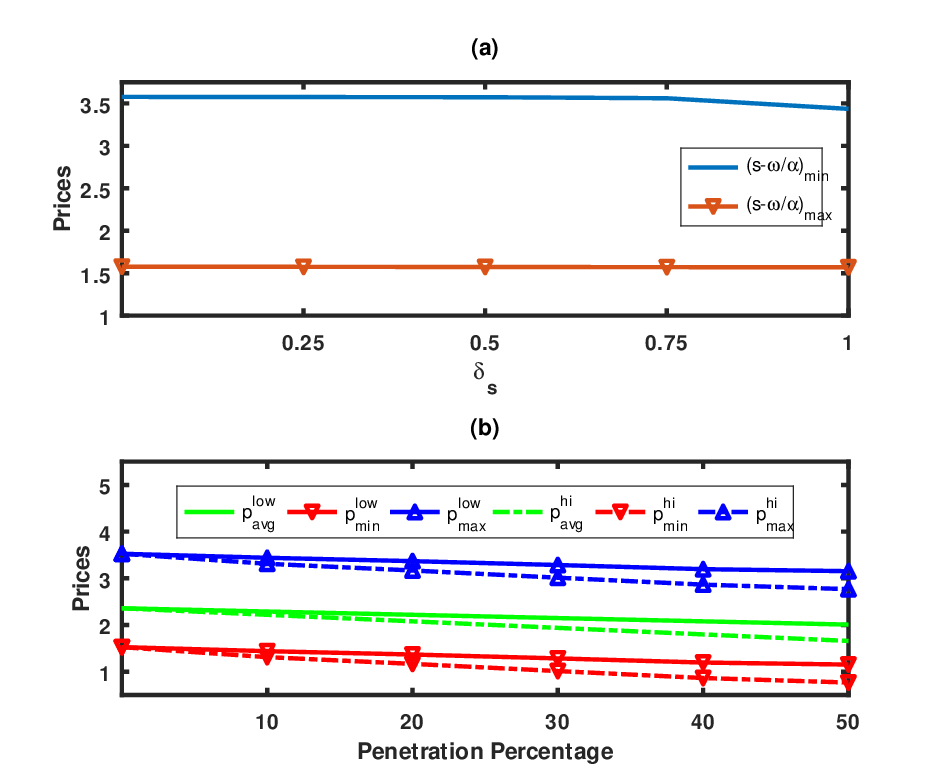}
    \vspace{-0.2in}
    \caption{Effects on prices for increased penetration at an individual level(a) or at community level(b)}
    \label{fig:penetration_ind}
    \vspace{-0.15in}
\end{figure}

In the second set of experiments, we also investigate the effects of change in penetration levels of renewable energy at an individual household level and when a fraction of prosumers in the community increase their capacities. We set up experiments where we increase the penetration level of a single household by a fixed amount $\delta_s$ first. We observe that price decreases as $\delta_s$ increases. However, the results vary according to the choice of household. We sort prosumers according to order of magnitude of ($s_i-\frac{\omega_i}{\alpha}$). ($s_i-\frac{\omega_i}{\alpha}$) represents the minimum surplus energy that the prosumer can dispose of/sell without sacrificing its utility. It is observed that for prosumers who are already contributing very high amounts of renewable energy to the grid, do not earn significant price cuts (price remains almost unchanged as $\delta_s$ is increased) and so they do not have any incentive to invest further in increasing penetration. On the contrary, prosumers who have low surplus to start with, get much larger price benefits (larger reductions in price) as they invest more in increasing penetration levels (increase $\delta_s$). Refer to Fig \ref{fig:penetration_ind}(a). However, increased penetration on an individual level does not affect prices of other prosumers much.     
 Fig \ref{fig:penetration_ind}(b) shows that as the fraction of prosumers within a community increase their renewable energies, the average price across the community decreases. A larger penetration level (higher $\delta_s$) results in a larger decrease in average prices. This conclusively shows that when a community invests in renewable energy, the entire community benefits from it in terms of prices.  

\section{Conclusion}
In this paper, we have considered the problem where the retailer can choose different prices for different users in a smart grid.  We propose a stylised model where the retailer's objective is to maximize revenue, minimize operation costs and also ensure high utility levels for the users whereas the user's objective is to maximize individual utilities. The pricing model is formulated as a Stackleberg game where the retailer and the users form a leader-follower pair. The retailer first sets a price and users respond by deciding their optimal amount of consumption. Since the original problem becomes non-convex, we propose 2 relaxed versions of the same problem which are convex and easier to solve. We  evaluate the optimality gaps across them in the asymptotic sense. We also explore the effects of price discrimination on retailer revenues and consumer utilities. We show that discrimination leads to a more equitable distribution of energy in the community, without appreciable change in net load. Retailer revenues are also found to increase with discrimination levels which makes such a pricing scheme lucrative to her. We extend our model to a scenario where users have in-house renewable energy generation capabilities and can feed energy back to the grid. A net-metering scenario is explored where the retail and the sell-back rates are the same and optimal prices are chosen by the retailer. Once again, we prove that discrimination leads to \textit{fairer pricing} where users are given price benefits when they contribute energy to the grid for similar willingness for demand. Price discrimination is also found to play an important role in incentivising more users to invest in renewable energy.  

\section{Future Work}
There are several interesting directions where our work can be extended. Throughout our analysis, we have assumed that consumption across different periods of the day are independent of each other. The characterization of prices when utilities of consumers are temporally correlated can be explored in the future. Another logical extension is exploring the pricing game in an incomplete information setup where there is uncertainty about some of the parameter values. In the scope of this paper, we have not considered storage as a viable option for a prosumer who generates energy in-house. Storage adds a new paradigm to the pricing problem in a net-metering scenario and is another interesting avenue where our work can be extended.   

\section{Acknowledgement}
This research did not receive any specific grant from funding agencies in the public, commercial, or not-for-profit sectors.

\section{Appendix}\label{sec:appendix}

\subsection{Proof for Theorem 1}
Let $p_1^{*}$ and $p_2^{*}$ be the optimal prices proposed by formulations 1 and 2 respectively. So, we seek to show that $p_2^{*} \geq p_1^{*}$ when $\eta = 0$.
\begin{proof}
Let $p^{*}$ be the unique global maximizer of the objective function in formulation 1 denoted by $f_1(p)$. $f_2(p)$ denotes the objective function in formulation 2. Let $min_{i=1(1)N}\{\omega_i\}$ be denoted by $\omega_L$. Now, there can be the following cases: a) $p^{*} \leq p_b$ b) $p_b < p^{*} \leq \omega_L$ c) $p^{*} > \omega_L$. 

For case a), because of the concavity of $f_1(\cdot)$, we must have $p_1^{*} = p_b$. Note that at $p_1^{*}$, $f_1(p_1^{*}) = f_2(p_1^{*})$. Let if possible, $p_2^{*} > p_1^{*}$. In that case, $f_2(p_2^{*}) = f_1(p_2^{*}) + \sum_{i}Min(0, \frac{\omega_i-p_2^{*}}{\alpha})$. Again due to concavity of $f_1(\cdot)$, $f_1(p_2^{*}) < f_1(p_1^{*})$. Also, we always have $\sum_{i}Min(0, \frac{\omega_i-p_2^{*}}{\alpha}) \leq 0$.  Therefore, $f_2(p_2^{*}) < f_1(p_1^{*}) = f_2(p_1^{*})$ which implies that $p_2^{*}$ cannot be optimal. So $p_2^{*} = p_1^{*} = p_b$.

For case b), we have $p_1^{*} = p^{*}$. Using a similar argument as in case a), it can be easily shown that $p_2^{*} = p_1^{*} = p^{*}$.

For case c), $p_1^{*} = \omega_L$. Let if possible, $p_2^{*} < p_1^{*}$. This implies, because of the concavity of $f_1(\cdot)$ and $p^{*} > p_1^{*}$, we must have $f_1(p_2^{*}) < f_1(p_1^{*})$. We have already argued that $f_2(p) \leq f_1(p)$ $\forall$ $p \in \mathbbm{R^{+}}$. Therefore, $f_2(p_2^{*}) < f_1(p_1^{*}) = f_2(p_1^{*})$ and $p_2^{*} < p_1^{*}$ cannot be optimal. So, we must have $p_2^{*} \geq p_1^{*}$. 

Combining all possible cases, it is clear that we will always have $p_2^{*} \geq p_1^{*}$. 
\end{proof}

\subsection{Proof for Theorem 2}
We seek to show that for any two consumers $i$ and $j$, if $\omega_i \geq \omega_j$, then in the optimal price vector, we must have $p_i \geq p_j$. Additionally, $\omega_i = \omega_j \implies p_i = p_j$. The proof has been done in two parts, in part 1 we prove for strict inequality and in part 2 we tackle the equality case. 
\begin{proof}
Let $(p_1, p_2,..p_i,..p_j,..p_n)$ be the optimal price vector obtained from the discriminatory pricing model. Let if possible, there exist a pair $i, j$ such that $\omega_i > \omega_j$, but $p_i < p_j$. Let $f_1(.)$ denote the retailer objective function in formulation 1.
If \textbf{p} is the optimal price vector, $f_1(\textbf{p}) \geq f_1(p')$ $\forall \hspace{3pt} p'\neq \textbf{p}$.\\
Now, let us consider a slightly modified price vector where prices $p_i$ and $p_j$ are interchanged. The new price vector is still a feasible solution (can be trivially verified). We will refer to this new price vector as $\textbf{q}$. 
\begin{equation*}
    \begin{split}
     f(\textbf{q}) - f(\textbf{p}) 
      & = \frac{e_1}{\alpha}(p_j(\omega_i-p_j)+p_i(\omega_j-p_i)-p_i(\omega_i-p_i)-p_j(\omega_j-p_j)) \\
      & = \frac{e_1}{\alpha}(p_j \omega_i + p_i \omega_j - p_i \omega_i - p_j \omega_j) = \frac{e_1}{\alpha}(\omega_i-\omega_j)(p_j-p_i) 
    \end{split}
\end{equation*}
When $\omega_i > \omega_j$, $f(\textbf{q})$ is strictly greater than $f(\textbf{p})$. This contradicts our initial assumption that $\textbf{p}$ is the optimal price vector. Hence, $p_i \geq p_j$ \hspace{3pt}$\forall\hspace{3pt}\omega_i > \omega_j$.

The objective in formulation 2 is given by $f_2(p) = f_1(p) + \sum_{k} Min(0, \frac{\omega_k-p_k}{\alpha})$ where $p$ is a n-dimensional price vector. To extend the result to formulation 2, we need to prove that $f_2(\textbf{q}) > f_2(\textbf{p})$ for the same choice of $\textbf{p}$ and $\textbf{q}$ as above. This essentially reduces to proving the following :
\begin{equation*}
    Min(0, \frac{\omega_i-p_j}{\alpha}) + Min(0, \frac{\omega_j-p_i}{\alpha}) \geq Min(0, \frac{\omega_i-p_i}{\alpha})+ Min(0, \frac{\omega_j-p_j}{\alpha}) 
\end{equation*}
when $\omega_i > \omega_j$ and $p_i < p_j$. Let $\Delta = Min(0, \frac{\omega_i-p_j}{\alpha}) + Min(0, \frac{\omega_j-p_i}{\alpha})-Min(0, \frac{\omega_i-p_i}{\alpha})- Min(0, \frac{\omega_j-p_j}{\alpha})$. $\omega_i$, $\omega_j$, $p_i$ and $p_j$ can be related in 24 ways. Because of the already assumed inequalities, there are 6 possible ways of arrangement. They are as follows :
\begin{itemize}
    \item $\omega_i > \omega_j > p_j > p_i$ : $\Delta = 0$
    \item $\omega_i > p_j > \omega_j > p_i$ : $\Delta = \frac{p_j-\omega_j}{\alpha} > 0 $
    \item $p_j > \omega_i > \omega_j > p_i$ : $\Delta = \frac{\omega_i-\omega_j}{\alpha} > 0$
    \item $\omega_i > p_j > p_i > \omega_j$ : $\Delta = \frac{p_j-p_i}{\alpha} > 0$
    \item $p_j > \omega_i > p_i > \omega_j$ : $\Delta = \frac{\omega_i-p_i}{\alpha} > 0$
    \item $p_j > p_i > \omega_i > \omega_j$ : $\Delta = 0$
\end{itemize}
Thus, we arrive at a contradiction and hence, our claim is valid, i.e., $p_i \geq p_j \forall$ $\omega_i > \omega_j$.
\end{proof}

\begin{proof}
For Part 2, we adopt a very similar approach as in Part 1. Let $\textbf{p}$ be the optimal price vector that maximizes $f_1(.)$. Thus, $f_1(\textbf{p}) \geq f_1(p')$ $\hspace{3pt}\forall\hspace{3pt}$ $p' \neq \textbf{p}$.
$\textbf{p}$ is given by $(p_1, p_2,..., p_i,.., p_j,...p_n)$ where $\omega_i = \omega_j$, but $p_i \neq p_j$. We now construct another price vector $\textbf{q}$ given by $(p_1, p_2,...\frac{(p_i+p_j)}{2},...\frac{(p_i+p_j)}{2},... p_n)$. Essentially, we have replaced prices $p_i$ and $p_j$ by $\frac{(p_i+p_j)}{2}$, the rest remain unchanged. $\textbf{q}$ is also a feasible solution (can be checked trivially). To prove by contradiction that our claim is correct, we need to show that $f_1(\textbf{q}) > f_1(\textbf{p})$. Let $\omega_i = \omega_j = \omega$. 
\begin{equation*}
    \begin{split}
        f_1(\textbf{p}) = e_1 (p_i\frac{(\omega-p_i)}{\alpha}+p_j\frac{(\omega-p_j)}{\alpha}) - \frac{e_2}{N} (K + \frac{\omega-p_i}{\alpha} + \frac{\omega-p_j}{\alpha})^{2} - e_3(\frac{p_i^{2}+p_j^{2}}{\alpha^{2}})  
    \end{split}
\end{equation*}
\begin{equation*}
    \begin{split}
        f_1(\textbf{q}) = e_1(p_i+p_j)\frac{(\omega-\frac{(p_i+p_j)}{2})}{\alpha} - \frac{e_2}{N}(K + \frac{2(\omega-\frac{(p_i+p_j)}{2})}{\alpha} )^{2} - 2 e_3(\frac{p_i+p_j}{2\alpha})^{2}  
    \end{split}
\end{equation*}
In order to prove the inequality, we do a term-by-term comparison. Observe that the second term in both $f_1(\textbf{p})$ and $f_1(\textbf{q})$ are the same and so they are not considered. We will use the following result for the proof :
\begin{equation*}
        AM > GM : p_i^{2} + p_j^{2} > 2 p_i p_j
        \implies 2(p_i^{2} + p_j^{2}) > (p_i + p_j)^{2} 
\end{equation*}
Therefore, for term 1,
\begin{equation*}
    \begin{split}
        (p_i+p_j)(\frac{\omega-\frac{(p_i+p_j)}{2}}{\alpha})
        = \omega(\frac{p_i+p_j}{\alpha}) - \frac{(p_i+p_j)^{2}}{2\alpha}
        \\> \omega(\frac{p_i+p_j}{\alpha}) - \frac{(p_i^{2}+p_j^{2})}{\alpha}
        = p_i(\frac{\omega-p_i}{\alpha})+p_j(\frac{\omega-p_j}{\alpha})
    \end{split}
\end{equation*}
And for term 3,
\begin{equation*}
    \frac{(p_i+p_j)^{2}}{2\alpha^{2}} < \frac{p_i^{2}+p_j^{2}}{\alpha^{2}}
\end{equation*}
Hence, $f_1(\textbf{q}) > f_1(\textbf{p})$ and the proof by contradiction is complete.\\
To extend the claim to formulation 2, we need to show additionally that:
\begin{equation*}
        2 Min(0, \frac{\omega - \frac{(p_i+p_j)}{2}}{\alpha}) \geq Min(0, \frac{\omega-p_i}{\alpha}) + Min(0, \frac{\omega-p_j}{\alpha})
\end{equation*}
Let $\delta = 2 Min(0, \frac{\omega - \frac{(p_i+p_j)}{2}}{\alpha}) - Min(0, \frac{\omega-p_i}{\alpha}) - Min(0, \frac{\omega-p_j}{\alpha})$. Without loss of generality, we can assume that $p_i > p_j$. That leaves us with 3 cases :
\begin{itemize}
    \item $\omega > p_i > p_j$ : $\delta = 0$
    \item $p_i > p_j > \omega$ : $\delta = 0$
    \item $p_i > \omega > p_j$ : This case has 2 sub-cases depending on whether $\omega > \frac{p_i+p_j}{2}$ or not. First we assume that $\omega$ is greater. Therefore, $\delta = \frac{p_i-\omega}{\alpha} > 0$. 
    When $\omega < \frac{p_i+p_j}{2}$, $\delta = \frac{\omega-p_j}{\alpha} > 0$.
\end{itemize}
Thus, the claim is proved to hold for formulation 2 as well. 
\end{proof}

\subsection{Asymptotic Analyses}
The entire asymptotic analysis is carried out in a non-discriminatory setting, so we have the same price $p$ for all consumers. We know that the elastic demand consumption $x_i$ is given by $x_i = Max(0, \frac{\omega_i-p}{\alpha})$. Now, we have two cases :  a) $p \leq \omega_{min}$ b) $p > \omega_{min}$. \\

\textit{Case a) :} When $p \leq \omega_{min}$, $x_i = \frac{\omega_i-p}{\alpha}$. Then we have the following :
\begin{equation*}
    E(\sum_{i}x_i) = \sum_{i}E(x_i) = \sum_{i} E(\frac{\omega_i - p}{\alpha}) = \frac{N(E(\omega) - p)}{\alpha}
\end{equation*}
\begin{equation*}
    \begin{split}
    E((\sum_{i}x_i)^2)  & = E^2(\sum_{i}x_i) + Var(\sum_{i}x_i) \\
    & = (\sum_{i}E(x_i))^2 + \sum_{i}Var(x_i) \hspace{10pt}(Cov(x_i, x_j) = 0 \hspace{3pt}\forall \hspace{3pt}i \neq j)\\
    & = \frac{N^{2}(E(\omega)-p)^{2}}{\alpha^2} + \frac{NVar(\omega)}{\alpha^{2}}
    \end{split}
\end{equation*}
\begin{equation*}
    E(\sum_{i} Min(0, \frac{\omega_i-p}{\alpha})) = E(0) = 0
\end{equation*}

\textit{Case b) :} When $p > \omega_{min}$, $x_i = Max(0, \frac{\omega_i-p}{\alpha})$. $\omega_i$'s are IID and drawn from $\mathbbm{U}[\omega_{min}, \omega_{max}]$. $f(\omega_i)$ is the PDF of $\omega_i$ and is given by $f(\omega_i) = \frac{1}{\omega_{max}-\omega_{min}}$. Then we have the following : 
\begin{equation*}
    \begin{split}
    E(x_i) & = E(Max(0, \frac{\omega_i-p}{\alpha})) = \int_{p}^{\omega_{max}}\frac{\omega_i-p}{\alpha}f(\omega_i)d\omega_i\\
    & = \frac{1}{\omega_{max}-\omega_{min}}\int_{p}^{\omega_{max}}\frac{\omega_i-p}{\alpha}d\omega_i = \frac{(\omega_{max}-p)^2}{2\alpha(\omega_{max}-\omega_{min})}
    \end{split}
\end{equation*}
\begin{equation*}
    \begin{split}
        E(x_i^2) & = \int_{p}^{\omega_{max}}(\frac{\omega_i-p}{\alpha})^{2}f(\omega_i)d\omega_i = \frac{(\omega_{max}-p)^3}{3\alpha^{2}(\omega_{max}-\omega_{min})}
    \end{split}
\end{equation*}
\begin{equation*}
     Var(x_i) = E(x_i^{2})-E^{2}(x_i) = \frac{(\omega_{max}-p)^{3}(\omega_{max}+3p-4\omega_{min})}{12\alpha^{2}(\omega_{max}-\omega_{min})^{2}}
\end{equation*}
\begin{equation*}
    \begin{split}
        E(Min(0, \frac{\omega_i-p}{\alpha})) & = \int_{\omega_{min}}^{p} \frac{\omega_i-p}{\alpha}f(\omega_i)d\omega_i\\
        & = \frac{1}{(\omega_{max}-\omega_{min})}\int_{\omega_{min}}^{p}\frac{\omega_i-p}{\alpha}d\omega_i = \frac{-(p-\omega_{min})^2}{2\alpha(\omega_{max}-\omega_{min})}
    \end{split}
\end{equation*}
Therefore, we have :
\begin{equation*}
    \begin{multlined}
        \lim_{N \to \infty}\frac{1}{N}E(g_0(p)) = \lim_{N \to \infty} e_1\cdot\frac{p(\omega_{max}-p)^2}{2\alpha(\omega_{max}-\omega_{min})} \\ - e_2\cdot\frac{(\omega_{max}-p)^{3}(\omega_{max}+3p-4\omega_{min})}{12N\alpha^{2}(\omega_{max}-\omega_{min})^{2}} - e_2\cdot\frac{(\omega_{max}-p)^{4}}{4\alpha^{2}(\omega_{max}-\omega_{min})^{2}} - e_3 \cdot \frac{p^2}{\alpha^2}\\
        = e_1\cdot\frac{p(\omega_{max}-p)^2}{2\alpha(\omega_{max}-\omega_{min})} - e_2\cdot\frac{(\omega_{max}-p)^{4}}{4\alpha^{2}(\omega_{max}-\omega_{min})^{2}} - e_3 \cdot \frac{p^2}{\alpha^2}
    \end{multlined}
\end{equation*}
\begin{equation*}
    \begin{multlined}
    \lim_{N \to \infty}\frac{1}{N}E(g_1(p)) = e_1\cdot\frac{p(E(\omega-p))}{\alpha} - e_2\cdot\frac{(E(\omega)-p)^{2}}{\alpha^2} - e_3\cdot\frac{p^2}{\alpha^2}
    \end{multlined}
\end{equation*}
\begin{equation*}
    \lim_{N \to \infty}\frac{1}{N}E(g_2(p)) = e_1\cdot\frac{p(E(\omega-p))}{\alpha} - e_2\cdot\frac{(E(\omega)-p)^{2}}{\alpha^2} - e_3\cdot\frac{p^2}{\alpha^2} - \frac{(p-\omega_{min})^2}{2\alpha(\omega_{max}-\omega_{min})}
\end{equation*}

\subsection{Uniqueness of solution for formulation 0(asymptotic)}
First, we need to check whether $g_0(p)$ is concave. $g_0(p)$ is given by $g_0(p) = K_1 p(\omega_{max}-p)^2 - K_2 (\omega_{max}-p)^4 - K_3 p^2$ and is biquadratic in $p$. In order for a biquadratic polynomial to be concave, its first order derivative (cubic polynomial) should have \textit{only one real root}. One possibility under which a cubic will always have only one real root, is if its derivative (quadratic) has no real roots. This leads us to the following condition under which $g_0(p)$ is always concave :  $(3K_1-2K_3)^2 < 24K_2\omega_{max}(K_3-K_1)$. Given our choice of $\alpha$, $\omega_{max}$ and $\omega_{min}$, this condition is very easily satisfied for most $(e_1, e_2, e_3)$ tuples where $e_1$, $e_2$ and $e_3$ $\in$ (0, 1].

\subsection{Uniqueness of solution for formulation 2(asymptotic)}
\begin{proof}
We have already argued that $g_1(p)$ and $g_2(p)$ are both concave in nature. Now we have, $p_2^{l} \in$ [$p_b$, $\omega_{min}$] and $p_2^{r} \in$ [$\omega_{min}$, $P$] such that $g_1(p_2^{l}) = g_2(p_2^{r})$. Clearly, $p_2^{r} \geq p_2^{l}$.

Now, $g_2(p_2^{r}) = g_1(p_2^{r}) - \frac{(p_2^{r}-\omega_{min})^2}{2\alpha(\omega_{max}-\omega_{min})}$. Therefore, $g_1(p_2^{l}) = g_2(p_2^{r})$ implies that $g_1(p_2^{r}) - g_1(p_2^{l}) = \frac{(p_2^{r}-\omega_{min})^2}{2\alpha(\omega_{max}-\omega_{min})} = K \geq 0$. (Condition 1)

There can be 2 possibilities : a) $p_2^{l}$ is equal to the maxima for $g_1$ which is obtained from the first order condition b) $p_2^{l}$ lies at one of the end-points of the interval [$p_b$, $\omega_{min}$] because the maxima lies outside the said interval.\\ \\
For case a), as $p_2^{l}$ is the unique global maximizer for $g_1$, $g_1(p_2^{l}) > g_1(q)$ for any $q$ in $\mathbbm{R^{+}}$. Therefore, for Condition 1 to be satisfied, $K = 0$ which implies that $p_2^{l} = p_2^{r} = \omega_{min}$. \\
For case b), the unique global maximizer for $g_1$ lies outside the range [$p_b$, $\omega_{min}$]. Since $p_2^{r} \geq p_2^{l}$ and $g_1(p_2^{r}) \geq g_1(p_2^{l})$, $p_2^{l} = \omega_{min}$ (because $g_1$ is concave). As $\omega_{min}$ also belongs to the interval [$\omega_{min}$, $P$], it is clear that we cannot find a price strictly to the left of the said interval such that $g_1(p_2^{l}) = g_2(p_2^{r})$. Therefore, the concavity of $g_2$ ensures uniqueness of the optimal solution of formulation 2 in the entire interval [$p_b$, $P$].     
\end{proof}

\bibliography{mybibfile}

\end{document}